\documentclass[12pt,a4paper,leqno]{article}

\usepackage[nottoc]{tocbibind}
\usepackage{latexsym}
\usepackage{amssymb,exscale}
\usepackage[centertags]{amsmath}
\usepackage{amsthm}
\usepackage{paralist}
\usepackage{footnpag}
\usepackage[all]{xy}

\usepackage[linktocpage]{hyperref}

\numberwithin{equation}{section}

\swapnumbers
\theoremstyle{definition}
\newtheorem{theorem}[equation]{Theorem}
\newtheorem{lemma}[equation]{Lemma}
\newtheorem{corollary}[equation]{Corollary}
\newtheorem{definition}[equation]{Definition}
\newtheorem{proposition}[equation]{Proposition}
\newtheorem{remark}[equation]{Remark}

\renewcommand{\phi}{\varphi}

\newcommand{\I}{{\rm i}}
\newcommand{\D}{\mathrm{d}}
\newcommand{\E}{\mathrm{e}}

\newcommand{\ti}{\tilde}

\renewcommand{\(}{\bigl(}
\renewcommand{\)}{\bigr)\vphantom{)}}

\newcommand{\ip}[2]{\langle#1,#2\rangle}

\newcommand{\imp}{$ \Longrightarrow $ }

\newcommand{\noise}{\text{\textup{noise}}}
\newcommand{\bin}{\text{\textup{bin}}}
\newcommand{\transfer}{\text{\textup{transfer}}}
\newcommand{\change}{\text{\textup{change}}}

\newcommand{\white}{\text{\textup{white}}}
\newcommand{\trivial}{\text{\textup{trivial}}}
\newcommand{\shift}{\text{\textup{shift}}}
\newcommand{\rotat}{\text{\textup{rotat}}}

\newcommand{\mes}{\operatorname{mes}}
\newcommand{\dist}{\operatorname{dist}}
\newcommand{\LocMin}{\operatorname{LocMin}}

\renewcommand{\Re}{\operatorname{Re}}
\renewcommand{\Im}{\operatorname{Im}}

\newcommand{\modO}{{\operatorname{mod}\,0}}

\newcommand{\One}{{1\hskip-2.5pt{\rm l}}}

\newcommand{\eps}{\varepsilon}
\newcommand{\si}{\sigma}
\newcommand{\ga}{\gamma}
\newcommand{\Ga}{\Gamma}
\newcommand{\om}{\omega}
\newcommand{\Om}{\Omega}

\newcommand{\al}{\alpha}
\newcommand{\be}{\beta}
\newcommand{\W}{\mathcal W}
\newcommand{\Ec}{\mathcal E}
\newcommand{\F}{\mathcal F}

\newcommand{\la}{\lambda}

\newcommand{\Ex}{\mathbb E\,}
\newcommand{\R}{\mathbb R}
\newcommand{\C}{\mathbb C}
\newcommand{\Z}{\mathbb Z}
\newcommand{\T}{\mathbb T}
\renewcommand{\Pr}[1]{\mathbb{P}\mskip1.5mu\(\mskip1.5mu#1\mskip1.5mu\)}
\newcommand{\cP}[2]{\mathbb{P}\mskip1.5mu\(\mskip1.5mu#1\mskip1.5mu
 \big|\mskip1.5mu#2\mskip1.5mu\)}

\newcommand{\sif}{$\sigma$\nobreakdash-field}

\newcommand{\invariant}[1]{$#1$\nobreakdash-\hspace{0pt}invariant}
\newcommand{\valued}[1]{$#1$\nobreakdash-\hspace{0pt}valued}
\newcommand{\measurable}[1]{$#1$\nobreakdash-\hspace{0pt}measurable}
\newcommand{\dimensional}[1]{$#1$\nobreakdash-\hspace{0pt}dimensional}

\setdefaultitem{\textasteriskcentered}{}{}{}
\makeatletter
\renewcommand*\l@section[2]{%
  \ifnum \c@tocdepth >\z@
    \addpenalty\@secpenalty
    \addvspace{0.25em \@plus\p@}%
    \setlength\@tempdima{2.5em}%
    \begingroup
      \parindent \z@ \rightskip \@pnumwidth
      \parfillskip -\@pnumwidth
      \leavevmode \bfseries
      \advance\leftskip\@tempdima
      \hskip -\leftskip
      #1\nobreak\hfil \nobreak\hb@xt@\@pnumwidth{\hss #2}\par
    \endgroup
  \fi}
\renewcommand*\numberline[1]{\hb@xt@\@tempdima{\hfil#1\hskip1em}}
\makeatother

\begin{document}

\title{On automorphisms of type $ I\! I $ Arveson systems\\
 (probabilistic approach)} 

\author{Boris Tsirelson}

\date{}
\maketitle

\begin{abstract}
A counterexample to the conjecture that the automorphisms of an arbitrary
Arveson system act transitively on its normalized units.
\end{abstract}

\setcounter{tocdepth}{1}
\tableofcontents

\pagebreak

\section*{Introduction}
\addcontentsline{toc}{section}{Introduction}
\smallskip

\parbox{6cm}{\small\textit{%
We do not know how to calculate the gauge group in this generality\dots}\\
\mbox{}\hfill W.~Arveson \cite[Sect.~2.8]{Ar}
}\hfill
\parbox{7.1cm}{\small\textit{%
At the moment, most important seems to us to answer the question whether the
automorphisms of an arbitrary product system act transitively on the
normalized units.}\\
\mbox{}\hfill V.~Liebscher \cite[Sect.~11]{Li}
}

\medskip

By an \emph{Arveson system} I mean a product system as defined by Arveson
\cite[3.1.1]{Ar}. Roughly, it consists of Hilbert spaces $ H_t $ (for $
0<t<\infty $) satisfying $ H_s \otimes H_t =  H_{s+t} $. Classical examples are
given by Fock spaces; these are \emph{type $ I $} systems, see \cite[3.3 and
Part 2]{Ar}. Their automorphisms are described explicitly, see
\cite[3.8.4]{Ar}. The group of automorphisms, called the \emph{gauge group} of
the Arveson system, for type $ I $ is basically the group of motions of the
\dimensional{N} Hilbert space. The parameter $ N \in \{0,1,2,\dots\} \cup
\{\infty\} $ is the so-called (numerical) index; accordingly, the system is
said to be of type $ I_0, I_1, I_2, \dots $ or $ I_\infty $. All Hilbert
spaces are complex (that is, over $ \C $).

Some Arveson systems contain no type $ I $ subsystems; these are \emph{type $
III $} systems, see \cite[Part 5]{Ar}. An Arveson system is of \emph{type $ II
$,} if it is not of type $ I $, but contains a type $ I $ subsystem. (See
\cite[6g and 10a]{Ts04} for examples.) In this case the greatest type $ I $
subsystem exists and will be called the \emph{classical part} of the type $ II
$ system. The latter is of type $ II_N $ where $ N $ is the index of its
classical part.

Little is known about the gauge group of a type $ II $ system and its natural
homomorphism into the gauge group of the classical part. In general, the
homomorphism is not one-to-one, and its range is a proper subgroup. The
corresponding subgroup of motions need not be transitive, which is the main
result of this work (Theorem \ref{1.8}); it answers a question asked by
Liebscher \cite[Notes 3.6, 5.8 and Sect. 11 (question 1)]{Li} and (implicitly)
Bhat \cite[Def. 8.2]{Bh01}; see also \cite{Ts04}, Question 9d3 and the
paragraph after it.

Elaborate constructions (especially, counterexamples) in a Hilbert space often
use a coordinate system (orthonormal basis). In other words, the sequence
space $ l_2 $ is used rather than an abstract Hilbert space. An Arveson system
consists of Hilbert spaces, but we cannot choose their bases without
sacrificing the given tensor product structure. Instead, we can choose maximal
commutative operator algebras, which leads to the probabilistic
approach. Especially, the white noise (or Brownian motion) will be used rather
than an abstract type $ I_1 $ Arveson system.

\section[Definitions, basic observations, and the result formulated]
 {\raggedright Definitions, basic observations, and the result formulated}
\label{sec:1}
I do not reproduce here the definition of an Arveson system \cite[3.1.1]{Ar},
since we only need the special case
\begin{equation}
  \label{eq:1.1}
  H_t = L_2 (\Om,\F_{0,t},P)
\end{equation}
corresponding to a noise.

\begin{definition}\label{1.2}
  A \emph{noise} consists of a probability space $
(\Om,\F,P) $, sub-\sif s $ \F_{s,t} \subset \F $ given for all $
s,t \in \R $, $ s<t $, and a measurable action $ (T_h)_h $ of $ \R $
on $ \Om $, having the following properties:
\begin{gather}
\F_{r,s} \otimes \F_{s,t} = \F_{r,t} \quad \text{whenever } r<s<t \,
 , \tag{a} \\
T_h \text{ sends } \F_{s,t} \text{ to } \F_{s+h,t+h} \quad
 \text{whenever } s<t \text{ and } h \in \R \, , \tag{b} \\
\F \text{ is generated by the union of all } \F_{s,t} \, . \tag{c}
\end{gather}
\end{definition}

See \cite[3d1]{Ts04} for details. As usual, all probability spaces are
standard, and everything is treated $ \modO $. Item (a) means that $ \F_{r,s}
$ and $ \F_{s,t} $ are (statistically) independent and generate $ \F_{r,t}
$. Invertible maps $ T_h : \Om \to \Om $ preserve the measure $ P $.

The white noise is a classical example; we denote it $ (\Om^\white,
\F^\white, P^\white) $, $ (\F_{s,t}^\white)_{s<t} $, $ (T_h^\white)_h $. It
is generated by the increments of the one-dimensional Brownian motion $
(B_t)_{-\infty<t<\infty} $, $ B_t : \Om \to \R $.

Given a noise, we construct Hilbert spaces $ H_t $ consisting of
\measurable{\F_{0,t}} complex-valued random variables, see \eqref{eq:1.1}. The
relation $ H_s \otimes H_t = H_{s+t} $, or rather a unitary operator $ H_s
\otimes H_t \to H_{s+t} $, emerges naturally,
\begin{multline*}
  H_{s+t} = L_2 (\F_{0,s+t}) = L_2 (\F_{0,s} \otimes \F_{s,s+t}) = \\
  = L_2 (\F_{0,s}) \otimes L_2(\F_{s,s+t}) = L_2 (\F_{0,s}) \otimes
  L_2(\F_{0,t}) = H_s \otimes H_t \, ;
\end{multline*}
the time shift $ T_s $ is used for turning $ \F_{s,s+t} $ to $ \F_{0,t}
$. Thus, $ (H_t)_{t>0} $ is an Arveson system. Especially, the white noise
leads to an Arveson system $ (H_t^\white)_{t>0} $ (of type $ I_1 $, as will be
explained).

For $ X \in H_s $, $ Y \in H_t $ the image of $ X \otimes Y $ in $ H_{s+t}
$ will be denoted simply $ XY $ (within this section).

We specialize the definition of a unit \cite[3.6.1]{Ar} to systems of the form
\eqref{eq:1.1}.

\begin{definition}
  A \emph{unit} (of the system \eqref{eq:1.1}) is a family $ (u_t)_{t>0} $ of
  non-zero vectors $ u_t \in H_t = L_2(\F_{0,t}) \subset L_2(\F) $ such that $
  t \mapsto u_t $ is a Borel measurable map $ (0,\infty) \to L_2(\F) $, and
\[
u_s u_t = u_{s+t} \quad \text{for all } s,t > 0 \, .
\]
\end{definition}

(In other words, the given unitary operator $ H_s \otimes H_t \to H_{s+t} $
maps $ u_s \otimes u_t $ to $ u_{s+t} $.) The unit is \emph{normalized,} if $
\| u_t \| = 1 $ for all $ t $. (In general, $ \| u_t \| = \exp(ct) $ for some
$ c \in \R $.)

Here is the general form of a unit in $ (H_t^\white)_t $:
\[
u_t = \exp ( z B_t + z_1 t ) \, ; \quad z,z_1 \in \C \, ;
\]
it is normalized iff $ (\Re z)^2 + \Re z_1 = 0 $. The units generate $
(H_t^\white)_t $ in the following sense: for every $ t > 0 $, $ H_t^\white $
is the closed linear span of vectors of the form $ (u_1)_{\frac t n}
(u_2)_{\frac t n} \dots (u_n)_{\frac t n} $, where $ u_1, \dots, u_n $ are
units, $ n=1,2,\dots $. Indeed, $ L_2 (\F_{0,t}) $ is spanned by random
variables of the form $ \exp \( \I \int_0^t f(s) \, \D B_s \) $ where $ f $
runs over step functions $ (0,t) \to \R $ constant on $ \( 0, \frac1n t \),
\dots, \( \frac{n-1}n t, t \) $.

We specialize two notions, `type $ I $' and `automorphism', to systems of the
form \eqref{eq:1.1}.

\begin{definition}
  A system of the form \eqref{eq:1.1} is of \emph{type} $ I $, if it is
  generated by its units.
\end{definition}

We see that $ (H_t^\white)_t $ is of type $ I $.

\begin{definition}\label{1.6}
  An \emph{automorphism} (of the system \eqref{eq:1.1}) is a family $
  (\Theta_t)_{t>0} $ of unitary operators $ \Theta_t : H_t \to H_t $ such that
  $ \Theta_{s+t} (XY) = (\Theta_s X) (\Theta_t Y) $ for all $ X \in H_s $, $ Y
  \in H_t $, $ s > 0 $, $ t > 0 $, and the function $ t \mapsto \ip{\Theta_t
  X_t}{Y_t} $ is Borel measurable whenever $ t \mapsto X_t $ and $ t \mapsto
  Y_t $ are Borel measurable maps $ (0,\infty) \to L_2(\F) $ such that $ X_t,
  Y_t \in L_2(\F_{0,t}) \subset L_2(\F) $.
\end{definition}

Basically, $ \Theta_s \otimes \Theta_t = \Theta_{s+t} $. The group $ G $ of all
automorphisms is called the \emph{gauge group.} Clearly, $ G $ acts on the set
of normalized units, $ (u_t)_t \mapsto (\Theta_t u_t)_t $.

Automorphisms $ \Theta_t = \Theta_t^{\trivial(\la)} = \E^{\I\la t} $ (for $
\la \in \R $), consisting of scalar operators, will be called \emph{trivial;}
these commute with all automorphisms, and are a one-parameter subgroup $
G^\trivial \subset G $. Normalized units $ (u_t)_t $ and $ (\E^{\I\la t}
u_t)_t $ will be called \emph{equilavent.} The factor group $ G / G^\trivial $
acts on the set of all equivalence classes of normalized units.

We turn to the gauge group $ G^\white $ of the classical system $
(H_t^\white)_t $. Equivalence classes of normalized units of $ (H_t^\white)_t
$ are parametrized by numbers $ z \in \C $, since each class contains exactly
one unit of the form
\[
u_t = \exp \( z B_t - (\Re z)^2 t \) \, .
\]
The scalar product corresponds to the distance:
\[
| \ip{ u_t^{(1)} }{ u_t^{(2)} } | = \exp \( -\tfrac12 |z_1-z_2|^2 t \)
\]
for $ u_t^{(k)} = \exp \( z_k B_t - (\Re z_k)^2 t \) $, $ k=1,2 $. The action
of $ G^\white / G^\trivial $ on equivalence classes boils down to its action
on $ \C $ by isometries. The orientation of $ \C $ is preserved, since
\[
\frac{ \ip{ u_t^{(1)} }{ u_t^{(2)} } \ip{ u_t^{(2)} }{ u_t^{(3)} } \ip{
    u_t^{(3)} }{ u_t^{(1)} } }{ | \ip{ u_t^{(1)} }{ u_t^{(2)} } \ip{ u_t^{(2)}
    }{ u_t^{(3)} } \ip{ u_t^{(3)} }{ u_t^{(1)} } | } =
\exp \( \I t S(z_1,z_2,z_3) \) \, ,
\]
where $ S(z_1,z_2,z_3) = \Im \( (z_2-z_1) \overline{ (z_3-z_1) } \) $ is twice
the \emph{signed} area of the triangle. So, $ G^\white / G^\trivial $ acts on
$ \C $ by motions (see \cite[3.8.4]{Ar}).

Shifts of $ \C $ along the imaginary axis, $ z \mapsto z + \I\la $ (for $ \la
\in \R $) emerge from automorphisms
\[
\Theta_t = \Theta_t^{\shift(\I\la)} = \exp (\I\la B_t) \, ;
\]
here the random variable $ \exp (\I\la B_t) \in L_\infty (\F_{0,t}^\white) $
is treated as the multiplication operator, $ X \mapsto X \exp (\I\la B_t) $
for $ X \in L_2(\F_{0,t}^\white) $.

Shifts of $ \C $ along the real axis, $ z \mapsto z + \la $ (for $ \la \in \R
$) emerge from less evident automorphisms 
\begin{equation}\label{1.star}
\Theta_t^{\shift(\la)} X = D_t^{1/2} \cdot ( X \circ \theta_t^\la ) \, ;
\end{equation}
here $ \theta_t^\la : C[0,t] \to C[0,t] $ is the drift transformation $
( \theta_t^\la b ) (s) = b(s) - 2\la s $ (for $ s \in [0,t] $), $ D_t $ is the
Radon-Nikodym derivative of the Wiener measure shifted by $ \theta_t^\la $
w.r.t.\ the Wiener measure itself,
\begin{equation}\label{1.65}
D_t = \exp ( 2\la B_t - 2\la^2 t ) \, ,
\end{equation}
and $ X \in L_2 (\F_{0,t}^\white) $ is treated as a function on $ C[0,t] $
(measurable w.r.t.\ the Wiener measure). Thus,
\[
( \Theta_t^{\shift(\la)} X ) (b) = \exp \( \la b(t) - \la^2 t \) X
(\theta_t^\la b) \, .
\]

By the way, these two one-parameter subgroups of $ G^\white $ satisfy Weyl
relations
\[
\Theta_t^{\shift(\la)} \Theta_t^{\shift(\I\mu)} = \E^{-2\I \la\mu t}
\Theta_t^{\shift(\I\mu)} \Theta_t^{\shift(\la)} \, ;
\]
that is, $ \Theta^{\shift(\la)} \Theta^{\shift(\I\mu)} =
\Theta^{\trivial(-2\la\mu)} \Theta^{\shift(\I\mu)} \Theta^{\shift(\la)} $.

Rotations of $ \C $ around the origin, $ z \mapsto \E^{\I\la} z $ (for $ \la
\in \R $) emerge from automorphisms $ \Theta^{\rotat(\la)} $. These will not
be used, but are briefly described anyway. They preserve Wiener chaos spaces $
H_n $,
\[
\Theta_t^{\rotat(\la)} X = \E^{\I n\la} X \quad \text{for } X \in H_n \cap
L_2(\F_{0,t}^\white) \, ;
\]
the $ n $-th chaos space $ H_n \subset L_2(\F^\white) $ consists of stochastic
integrals
\[
X = \idotsint_{-\infty<s_1<\dots<s_n<\infty} f(s_1,\dots,s_n) \, \D B_{s_1}
\dots \D B_{s_n}
\]
where $ f \in L_2(\R^n) $ (or rather, the relevant part of $ \R^n $). One may
say that $ \Theta^{\rotat(\la)} $ just multiplies each $ \D B_s $ by $
\E^{\I\la} $.

Combining shifts and rotations we get all motions of $ \C $. Accordingly, all
automorphisms of $ (H_t^\white)_t $ are combinations of $
\Theta^{\shift(\I\la)} $, $ \Theta^{\shift(\la)} $, $ \Theta^{\rotat(\la)} $
and $ \Theta^{\trivial(\la)} $. More generally, the \dimensional{N} Brownian
motion leads to the (unique up to isomorphism) Arveson system of type $ I_N
$ and motions of $ \C^N $. We need $ N=1 $ only; $ (H_t^\white)_t $ is the
Arveson system of type $ I_1 $.

Some noises are constructed as extensions of the white noise,
\begin{equation}\label{eq:1.65}
\F_{s,t} \supset \F_{s,t}^\white
\end{equation}
(also $ T_h $ conforms to $ T_h^\white $). More exactly, it means that $ B_t
\in L_2(\F) $ are given such that $ B_t - B_s $ is \measurable{\F_{s,t}} for $
-\infty < s < t < \infty $, and $ B_t - B_s \sim N(0,t-s) $ (that is, the
random variable $ (t-s)^{-1/2} (B_t - B_s) $ has the standard normal
distribution), and $ B_0 = 0 $, and $ (B_t - B_s) \circ T_h = B_{t+h} -
B_{s+h} $. Such $ (B_t)_t $ may be called a Brownian motion adapted to the
given noise. Then, of course, by $ \F_{s,t}^\white $ we mean the sub-\sif\
generated by $ B_u - B_s $ for all $ u \in (s,t) $. The Arveson system $
(H_t)_t $, $ H_t = L_2(\F_{0,t}) $, is an extension of the type $ I_1 $ system
$ (H_t^\white)_t $, $ H_t^\white = L_2(\F_{0,t}^\white) $,
\begin{equation}\label{eq:1.67}
H_t \supset H_t^\white \, .
\end{equation}
All units of $ (H_t^\white)_t $ are also units of $ (H_t)_t $. It may happen
that $ (H_t)_t $ admits no other units even though $ \F_{s,t} \ne
\F_{s,t}^\white $, $ H_t \ne H_t^\white $. Then $ (H_t)_t $ is of type $ II $
(units generate a nontrivial, proper subsystem), namely, of type $ II_1 $; $
(H_t^\white)_t $ is the classical part of $ (H_t)_t $, and the white noise is
the classical part of the given noise. The automorphisms $
\Theta^\trivial(\la) $ and $ \Theta^{\shift(\I\la)} $ for $ \la \in \R $ can
be extended naturally from the classical part to the whole system (which does
not exclude other possible extensions). For $ \Theta^{\shift(\la)} $ and $
\Theta^{\rotat(\la)} $ we have no evident extension. Moreover, these
automorphisms need not have any extensions, as will be proved.

Two examples found by Warren \cite{Wa1}, \cite{Wa} are `the noise of
splitting' and 'the noise of stickiness'; see also \cite{Wat0} and
\cite[Sect.~2]{Ts04}. For the noise of splitting the gauge group restricted to
the classical part covers all shifts of $ \C $ (but only trivial rotations
\cite{Ts06}), thus, it acts transitively on $ \C $, therefore, on normalized
units as well.

A new (third) example is introduced in Sect.~\ref{sec:7} for proving the
main result formulated as follows.

\begin{theorem}\label{1.8}
There exists an Arveson system of type $ II_1 $ such that the action of the
group of automorphisms on the set of normalized units is not transitive.
\end{theorem}

The proof is given in Sect.~\ref{sec:8}, after the formulation of
Prop.~\ref{10.1}.

A weaker result, obtained by different methods, was reported \cite{Ma}.

The first version \cite{I} of this paper have raised some doubts
\cite[p.~6]{BLS}. Hopefully they will be dispelled by the present version.

First of all, in Sect.~\ref{sec:1-5} we reformulate the problem as a problem
of isomorphism. Isomorphism of some models simpler than Arveson systems are
investigated in Sections \ref{sec:2}--\ref{sec:6}. In Sect.~\ref{sec:8} we
reduce the problem for Arveson systems to the problem for the simpler
models. In combination with the new noise of Sect.~\ref{sec:7} it proves
Theorem \ref{1.8}.

\section[Extensions of automorphisms and isomorphisms of extensions]
 {\raggedright Extensions of automorphisms and isomorphisms of extensions}
\label{sec:1-5}
Assume that a given noise $ \( (\Om,\F,P), (\F_{s,t}), (T_h) \) $ is an
extension of the white noise (see \eqref{eq:1.65} and the explanation after
it) generated by a given Brownian motion $ (B_t)_t $ adapted to the given
noise. Assume that another noise $ \( (\Om',\F',P'), (\F'_{s,t}), (T'_h) \) $
is also an extension of the white noise, according to a given adapted Brownian
motion $ (B'_t)_t $. On the level of Arveson systems we have two extensions of
the type $ I_1 $ system:
\[
H_t \supset H_t^\white \, , \quad H'_t \supset H_t^{\prime\>\white} \, ;
\]
here $ H_t = L_2 (\Om,\F_{0,t},P) $, $ H_t^\white = L_2
(\Om,\F_{0,t}^\white,P) $, $ \F_{0,t}^\white $ being generated by the
restriction of $ B $ to $ [0,t] $. ($ H'_t $ and $ H_t^{\prime\>\white} $ are
defined similarly.)

An isomorphism between the two Arveson systems $ (H_t)_t $, $ (H'_t)_t $ is
defined similarly to \ref{1.6} ($ \Theta_t : H_t \to H'_t $, $ \Theta_{s+t} =
\Theta_s \otimes \Theta_t $, and the Borel measurability). If it exists, it is
non-unique. In contrast, the subsystems $ (H_t^\white)_t $ and $
(H_t^{\prime\>\white})_t $ are \emph{naturally} isomorphic:
\[
\Theta_t^\transfer \( X (B|_{[0,t]}) \) = X (B'|_{[0,t]}) \quad \text{for all
} X \, ;
\]
here $ B|_{[0,t]} $ is treated as a \valued{C[0,t]} random variable on $ \Om
$, distributed $ \W_t $ (the Wiener measure); similarly, $ B'|_{[0,t]} $ is
a \valued{C[0,t]} random variable on $ \Om' $, distributed $ \W_t $; and $ X
$ runs over $ L_2 (C[0,t], \W_t) $.

We define an \emph{isomorphism between extensions} as an isomorphism $
(\Theta_t)_t $ between Arveson systems that extends $ \Theta^\transfer $, that
is,
\[
\Theta_t |_{H_t^\white} = \Theta_t^\transfer \quad \text{for all } t \, .
\]

Adding a drift to the Brownian motion $ (B_t)_t $ we get a random process $
(B_t + \la t )_t $ locally equivalent, but globally singular to the Brownian
motion. In terms of noises this idea may be formalized as follows.

Let $ (\Om,\F,\ti P) $ be a probability space, $ \F_{s,t} \subset \F $
sub-\sif s, and $ (T_h)_h $ a measurable action of $ \R $ on $ \Om $,
satisfying Conditions (b) and (c) of Def.~\ref{1.2} (but not (a)).
Let $ P, P' $ be \invariant{(T_h)} probability measures on $ (\Om,\F) $
such that $ P + P' = 2 \ti P $, and \ref{1.2}(a) holds for each of the two
measures $ P, P' $. Then we have two noises $ \( (\Om,\F,P), (\F_{s,t}),
(T_h) \) $, $ \( (\Om,\F,P'), (\F_{s,t}), (T_h) \) $. Assume also that the
restrictions $ P |_{\F_{s,t}} $ and $ P' |_{\F_{s,t}} $ are equivalent (that
is, mutually absolutely continuous) whenever $ s<t $. This relation between
two noises may be called a change of measure. The corresponding Arveson
systems are naturally isomorphic (via multiplication by the Radon-Nikodym
derivative):
\[
\Theta_t^\change : H_t \to H'_t \, , \quad \Theta_t^\change \psi = D_t^{-1/2}
\psi \, , \quad D_t = \frac{ \D P'|_{\F_{0,t}} }{ \D P|_{\F_{0,t}} } \, .
\]

We are especially interested in a change of measure such that (recall
\eqref{1.65})
\[
D_t = \exp \( 2 \la B_t - 2 \la^2 t \) \quad \text{for } t \in (0,\infty) \, ,
\]
where $ (B_t)_t $ is a Brownian motion adapted to the first noise, and $ \la
\in \R $ a given number. In this case $ (B_t - 2 \la t)_t $ is a Brownian
motion adapted to the second noise. We take $ B'_t = B_t - 2\la t $ and get
two extensions of the white noise. In such a situation we say that the second
extension results from the first one by the drift $ 2\la $, denote $
\Theta_t^\change $ by $ \Theta_t^{\change(\la)} $ and $ \Theta_t^\transfer $
by $ \Theta_t^{\transfer(\la)} $.

Note that $ \Theta_t^{\transfer(\la)} \( X (B|_{[0,t]}) \) = ( X \circ
\theta_t^\la ) (B|_{[0,t]}) $ for $ X \in L_2 (C[0,t], \W_t) $; as before, $
\theta_t^\la : C[0,t] \to C[0,t] $ is the drift transformation, $ (
\theta_t^\la b ) (s) = b(s) - 2\la s $ for $ s \in [0,t] $, it sends the
measure $ D_t \cdot \W_t $ to $ \W_t $.

The isomorphism $ \Theta^{\change(\la)} $ between the two Arveson systems $
(H_t)_t $, $ (H'_t)_t $ is not an isomorphism of extensions (unless $ \la=0
$), since its restiction to $ (H_t^\white)_t $ is not equal to $
\Theta^{\transfer(\la)} $. Instead, by the lemma below, they are related via
the automorphism $ \Theta^{\shift(\la)} $ of $ (H_t^\white)_t $ introduced in
Sect.~\ref{sec:1}.

\begin{lemma}
\[
\Theta^{\change(\la)} \Theta^{\shift(\la)} = \Theta^{\transfer(\la)} \, ,
\]
that is,
\[
\Theta_t^{\change(\la)} \Theta_t^{\shift(\la)} \psi =
\Theta_t^{\transfer(\la)} \psi
\]
for all $ \psi \in H_t^\white $ and all $ t \in (0,\infty) $.
\end{lemma}

\begin{proof}
We take $ X \in L_2 (C[0,t], \W_t) $ such that $ \psi = X \( B|_{[0,t]} \)
$, then
\begin{multline*}
\Theta_t^{\change(\la)} \Theta_t^{\shift(\la)} \psi = D_t^{-1/2} \cdot
 D_t^{1/2} \cdot ( X \circ \theta_t^\la ) ( B|_{[0,t]} ) = \\
= ( X \circ \theta_t^\la ) ( B|_{[0,t]} ) = \Theta_t^{\transfer(\la)} \psi \,
.
\end{multline*}
\end{proof}

The situation is shown on the diagram
\[
\xymatrix@C+1cm{
 (H_t)_t \ar@{.>}[r]_{\Theta} \ar@{.>}@/^1.5pc/[rr]^{\Theta'} & (H_t)_t
 \ar[r]_{\Theta^{\change(\la)}} & (H'_t)_t
\\
 (H_t^\white)_t \ar@{^{(}->}[u] \ar[r]^{\Theta^{\shift(\la)}}
 \ar@/_1.5pc/[rr]_{\Theta^{\transfer(\la)}} & (H_t^\white)_t
 \ar@{^{(}->}[u] \ar[r]^{\Theta^{\change(\la)}} & (H_t^{\prime\>\white})_t
 \ar@{^{(}->}[u]
}
\]
and we see that the following conditions are equivalent:
\begin{compactitem}
\item there exists an automorphism $ \Theta $ of $ (H_t)_t $ that extends $
 \Theta^{\shift(\la)} $;
\item there exists an isomorphism $ \Theta' $ between $ (H_t)_t $ and $
 (H'_t)_t $ that extends $ \Theta^{\transfer(\la)} $.
\end{compactitem}
In other words, $ \Theta^{\shift(\la)} $ can be extended to $ (H_t)_t $ if and
only if the two extensions of the type $ I_1 $ system are isomorphic.

\begin{corollary}\label{cor2}
In order to prove Theorem \ref{1.8} it is sufficient to construct a noise,
extending the white noise, such that for every $ \la \in \R \setminus \{0\} $
the extension obtained by the drift $ \la $ is non-isomorphic to the original
extension on the level of Arveson systems (that is, the corresponding
extensions of the type $ I_1 $ Arveson system are non-isomorphic).
\end{corollary}

\begin{proof}
In the group of all motions of the complex plane we consider the subgroup $ G
$ of motions that correspond to automorphisms of $ (H_t^\white)_t $ extendable
to $ (H_t)_t $. Real shifts $ z \mapsto z+\la $ (for $ \la \in \R \setminus
\{0\} $) do not belong to $ G $, as explained above. Imaginary shifts $ z
\mapsto z+\I\la $ (for $ \la \in \R $) belong to $ G $, since the operators $
\Theta_t^{\shift(\I\la)} $ of multiplication by $ \exp(\I\la B_t) $ act
naturally on $ H_t $. It follows that $ G $ contains no rotations (except for
the rotation by $ \pi $) and therefore is not transitive.
\end{proof}

Thus, we need a \emph{drift sensitive} extension. Such extension is
constructed in Sect.~\ref{sec:7} and its drift sensitivity is proved in
Sect.~\ref{sec:8}.

\section[Toy models: Hilbert spaces]
 {\raggedright Toy models: Hilbert spaces}
\label{sec:2}
Definitions and statements of Sections \ref{sec:2} and \ref{sec:3} will not be
used formally, but probably help to understand the idea.

The phenomenon of a non-extendable isomorphism (as well as nonisomorphic
extensions) is demonstrated in this section by a toy model, --- a kind of
product system of Hilbert spaces, simpler than Arveson system.

\begin{definition}
A \emph{toy product system of Hilbert spaces} is a triple $ ( H_1, H_\infty,
U) $, where $ H_1, H_\infty $ are Hilbert spaces (over $ \C $, separable), and
$ U : H_1 \otimes H_\infty \to H_\infty $ is a unitary operator.
\end{definition}

We treat it as a kind of product system, since
\[
H_\infty \sim H_1 \otimes H_\infty \sim H_1 \otimes H_1 \otimes H_\infty \sim
\dots
\]
where `$ \sim $' means: may be identified naturally (using $ U $).

An evident example: $ H_\infty = (H_1,\psi_1)^{\otimes\infty} $ is the
infinite tensor product of (an infinite sequence of) copies of $ H_1 $
relatively to (the copies of) a given vector $ \psi_1 \in H_1 $, $ \| \psi_1
\| = 1 $. The equation $ U(\psi\otimes\xi) = \xi $ has exactly one solution: $
\psi = \psi_1 $, $ \xi = \psi_1^{\otimes\infty} $.

An uninteresting modification: $ H_\infty = (H_1,\psi_1)^{\otimes\infty}
\otimes H_0 $ for some Hilbert space $ H_0 $.

A more interesting example: $ H_\infty = (H_1,\psi_1)^{\otimes\infty} \oplus
(H_1,\psi_2)^{\otimes\infty} $ is the direct sum of two such infinite tensor
products, one relative to $ \psi_1 $, the other relative to another vector $
\psi_2 \in H_1 $, $ \| \psi_2 \| = 1 $, $ \psi_2 \ne \psi_1 $. The equation $
U(\psi\otimes\xi) = \xi $ has exactly two solutions: $ \psi = \psi_1 $, $ \xi
= \psi_1^{\otimes\infty} $ and $ \psi = \psi_2 $, $ \xi =
\psi_2^{\otimes\infty} $.

\begin{definition}
Let $ ( H_1, H_\infty, U) $ and $ ( H'_1, H'_\infty, U') $ be toy product
systems of Hilbert spaces. An \emph{isomorphism} between them is a pair $
\Theta = (\Theta_1,\Theta_\infty) $ of unitary operators $ \Theta_1 : H_1 \to
H'_1 $, $ \Theta_\infty : H_\infty \to H'_\infty $ such that the diagram
\[
\xymatrix@C+1cm{
 H_1 \otimes H_\infty \ar[d]^{\Theta_1 \otimes \Theta_\infty} \ar[r]^{U} &
  H_\infty \ar[d]^{\Theta_\infty}
\\
 H'_1 \otimes H'_\infty \ar[r]^{U'} & H'_\infty
}
\]
is commutative.
\end{definition}

Thus,
\[
\Theta_\infty \sim \Theta_1 \otimes \Theta_\infty \sim \Theta_1 \otimes
\Theta_1 \otimes \Theta_\infty \sim
\dots
\]

A unitary operator $ \Theta_1 : H_1 \to H_1 $ leads to an automorphism of $
(H_1,\psi_1)^{\otimes\infty} $ (that is, of the corresponding toy product
system) if and only if $ \Theta_1 \psi_1 = \psi_1 $. Similarly, $ \Theta_1 $
leads to an automorphism of $ (H_1,\psi_1)^{\otimes\infty} \oplus
(H_1,\psi_2)^{\otimes\infty} $ if and only if either $ \Theta_1 \psi_1 =
\psi_1 $ and $ \Theta_1 \psi_2 = \psi_2 $, or $ \Theta_1 \psi_1 = \psi_2 $ and
$ \Theta_1 \psi_2 = \psi_1 $.

Taking $ \Theta_1 $ such that $ \Theta_1 \psi_1 = \psi_1 $ but $ \Theta_1
\psi_2 \ne \psi_2 $ we get an automorphism of $ (H_1,\psi_1)^{\otimes\infty} $
that cannot be extended to an automorphism of $ (H_1,\psi_1)^{\otimes\infty}
\oplus (H_1,\psi_2)^{\otimes\infty} $.

Similarly to Sect.~\ref{sec:1-5} we may turn from extensions of automorphisms
to isomorphisms of extensions. The system $ (H_1,\psi_1)^{\otimes\infty}
\oplus (H_1,\psi_2)^{\otimes\infty} $ is an extension of $
(H_1,\psi_1)^{\otimes\infty} $ (in the evident sense). Another vector $
\psi'_2 $ leads to another extension of $ (H_1,\psi_1)^{\otimes\infty} $. We
define an isomorphism between the two extensions as an isomorphism $
(\Theta_1,\Theta_\infty) $ between the toy product systems $
(H_1,\psi_1)^{\otimes\infty} \oplus (H_1,\psi_2)^{\otimes\infty} $ and $
(H_1,\psi_1)^{\otimes\infty} \oplus (H_1,\psi'_2)^{\otimes\infty} $ whose
restriction to $ (H_1,\psi_1)^{\otimes\infty} $ is trivial (the identity).
\[
\xymatrix{
(H_1,\psi_1)^{\otimes\infty} \oplus (H_1,\psi_2)^{\otimes\infty}
 \ar@{<->}[rr]^{(\Theta_1,\Theta_\infty)} & &
 (H_1,\psi_1)^{\otimes\infty} \oplus (H_1,\psi'_2)^{\otimes\infty}
\\
 & (H_1,\psi_1)^{\otimes\infty} \ar@{_{(}->}[ul] \ar@{^{(}->}[ur]
}
\]
Clearly, $ \Theta_1 $ must be trivial; therefore $ \psi'_2 $ must be equal to
$ \psi_2 $. Otherwise the two extensions are nonisomorphic.

\section[Toy models: probability spaces]
 {\raggedright Toy models: probability spaces}
\label{sec:3}
\begin{definition}
A \emph{toy product system of probability spaces} is a triple $ ( \Om_1,
\Om_\infty,
\linebreak[0]
\al) $, where $ \Om_1, \Om_\infty $ are probability spaces (standard), and $
\al : \Om_1 \times \Om_\infty \to \Om_\infty $ is an isomorphism $ \modO $
(that is, an invertible measure preserving map).
\end{definition}

Every toy product system of probability spaces $ ( \Om_1, \Om_\infty, \al) $
leads to a toy product system of Hilbert spaces $ ( H_1, H_\infty, U) $ as
follows:
\begin{gather*}
H_1 = L_2 (\Om_1) \, ; \quad H_\infty = L_2 (\Om_\infty) \, ; \\
(U\psi)(\cdot) = \psi ( \al^{-1} (\cdot) ) \, .
\end{gather*}
Here we use the canonical identification
\[
L_2 (\Om_1) \otimes L_2 (\Om_\infty) = L_2 ( \Om_1 \times \Om_\infty )
\]
and treat a vector $ \psi \in H_1 \otimes H_\infty $ as an element of $ L_2 (
\Om_1 \times \Om_\infty ) $.

An evident example: $ \Om_\infty = \Om_1^\infty $ is the product of an
infinite sequence of copies of $ \Om_1 $. It leads to $ H_\infty =
(H_1,\One)^{\otimes\infty} $ where $ H_1 = L_2 (\Om_1) $ and $ \One \in
L_2(\Om_1) $ is the constant function, $ \One(\cdot)=1 $.

An uninteresting modification: $ \Om_\infty = \Om_1^\infty \times \Om_0 $ for
some probability space $ \Om_0 $. It leads to $ H_\infty =
(H_1,\One)^{\otimes\infty} \otimes H_0 $, $ H_0 = L_2(\Om_0) $.

Here is a more interesting example. Let $ X_1 : \Om_1 \to \{-1,+1\} $ be a
measurable function (not a constant). We define $ \Om_\infty $ as the set of
all double sequences $ \( \begin{smallmatrix} \om_1, & \om_2, & \dots \\ s_1,
& s_2, & \dots \end{smallmatrix} \) $ such that $ \om_k \in \Om_1 $, $ s_k \in
\{-1,+1\} $ and $ s_k = s_{k+1} X_1 (\om_k) $ for all $ k $. Sequences $
(\om_1,\om_2,\dots) \in \Om_1^\infty $ are endowed with the product
measure. The conditional distribution of the sequence $ (s_1,s_2,\dots) $,
given $ (\om_1,\om_2,\dots) $, must be concentrated on the two sequences
obeying the relation $ s_k = s_{k+1} X_1 (\om_k) $. We give to these two
sequences equal conditional probabilities, $ 0.5 $ to each. Thus, $ \Om_\infty
$ is endowed with a probability measure. The map $ \al : \Om_1 \times
\Om_\infty \to \Om_\infty $ is defined by
\[
\al \bigg( \om_1, \begin{pmatrix} \om_2, & \om_3, & \dots \\ s_2, & s_3, &
\dots \end{pmatrix} \bigg) = \begin{pmatrix} \om_1, & \om_2, & \om_3, & \dots
\\ s_2 X_1 (\om_1), & s_2, & s_3, & \dots \end{pmatrix} \, .
\]
Clearly, $ \al $ is measure preserving.

This system $ (\Om_1,\Om_\infty,\al) $ leads to a system $ (H_1,H_\infty,U) $
of the form $ (H_1,\psi_1)^{\otimes\infty} \oplus (H_1,\psi_2)^{\otimes\infty}
$ (up to isomorphism), as explained below. We have
\begin{gather*}
H_1 = L_2 (\Om_1) \, , \quad H_\infty = L_2 (\Om_\infty) \, , \\
(U\psi) \begin{pmatrix} \om_1, & \om_2, & \om_3, & \dots \\ s_1, & s_2, & s_3,
  & \dots \end{pmatrix} = \psi \bigg( \om_1, \begin{pmatrix} \om_2, & \om_3, &
  \dots \\ s_2, & s_3, & \dots \end{pmatrix} \bigg) \, .
\end{gather*}
The equation $ U(\psi\otimes\xi) = \xi $ becomes
\[
\psi(\om_1) \xi \begin{pmatrix} \om_2, & \om_3, & \dots \\ s_2, & s_3, & \dots
\end{pmatrix} = \xi \begin{pmatrix} \om_1, & \om_2, & \om_3, & \dots \\ s_1, &
s_2, & s_3, & \dots \end{pmatrix} \, .
\]
One solution is evident: $ \psi = \One_{\Om_1} $, $ \xi = \One_{\Om_\infty}
$. A less evident solution is, $ \psi = X_1 $, $ \xi = S_1 $, where $ S_1 $ is
defined by $ S_1 \( \begin{smallmatrix} \om_1, & \om_2, & \dots \\ s_1, &
s_2, & \dots \end{smallmatrix} \) = s_1 $. (The equation is satisfied due to
the relation $ X_1(\om_1) s_2 = s_1 $.) We consider the system $
(H'_1,H'_\infty,U') $ where $ H'_1 = H_1 = L_2(\Om_1) $, $ H'_\infty = (H'_1,
\One_{\Om_1})^{\otimes\infty} \oplus (H'_1,X_1)^{\otimes\infty} $ ($ U' $
being defined naturally) and construct an isomorphism $
(\Theta_1,\Theta_\infty) $ between $ (H_1,H_\infty,U) $ and $
(H'_1,H'_\infty,U') $ such that $ \Theta_\infty \One_{\Om_\infty} =
\One_{\Om_1}^{\otimes\infty} $, $ \Theta_\infty S_1 = X_1^{\otimes\infty}
$. To this end we consider an arbitrary $ n $ and $ \xi \in L_2(\Om_1^n) =
H_1^{\otimes n} $, define $ \phi, \psi \in L_2(\Om_\infty) $ by
\begin{align*}
\phi \begin{pmatrix} \om_1, & \om_2, & \dots \\ s_1, & s_2, & \dots
 \end{pmatrix} &= \xi (\om_1,\dots,\om_n) \, , \\
\psi \begin{pmatrix} \om_1, & \om_2, & \dots \\ s_1, & s_2, & \dots
 \end{pmatrix} &= s_{n+1} \xi (\om_1,\dots,\om_n)
\end{align*}
and, using the relation (or rather, the natural isomorphism) $ H'_\infty =
(H'_1)^{\otimes n} \otimes H'_\infty $, we let
\[
\Theta_\infty \phi = \xi \otimes \One_{\Om_1}^{\otimes\infty} \, , \quad
\Theta_\infty \psi = \xi \otimes X_1^{\otimes\infty} \, ,
\]
thus defining a unitary $ \Theta_\infty : H_\infty \to H'_\infty $. (Further
details are left to the reader.)

A more general construction is introduced in Sect.~\ref{sec:4}.

\section[Binary extensions: probability spaces]
 {\raggedright Binary extensions: probability spaces}
\label{sec:4}
\begin{definition}
(a) An \emph{extension} of a probability space $ \Om $ consists of another
probability space $ \ti\Om $ and a measure preserving map $ \ga : \ti\Om \to
\Om $.

(b) Two extensions $ (\ti\Om,\ga) $ and $ (\ti\Om',\ga') $ of a probability
space $ \Om $ are \emph{isomorphic,} if there exists an invertible ($ \modO $)
measure preserving map $ \theta : \ti\Om \to \ti\Om' $ such that the diagram
\[
\xymatrix{
 \ti\Om \ar[dr]_{\ga} \ar[rr]^{\theta} & & \ti\Om' \ar[dl]^{\ga'}
\\
 & \Om
}
\]
is commutative. (Such $ \theta $ will be called an \emph{isomorphism of
extensions}.)

(c) An extension of a probability space $ \Om $ is \emph{binary,} if it is
isomorphic to $ (\Om\times\Om_\pm,\ga) $, where $ \Om_\pm = \{-1,+1\} $
consists of two equiprobable atoms, and $ \ga : \Om\times\Om_\pm \to \Om $ is
the projection, $ (\om,s) \mapsto \om $.
\end{definition}

By a well-known theorem of V.~Rokhlin, an extension is binary if and only if
conditional measures consist of two atoms of probability $ 0.5 $. However,
this fact will not be used.

Interchanging the two atoms we get an involution on $ \ti\Om $. Denoting it by
$ \ti\om \mapsto -\ti\om $ we have
\[
-\ti\om \ne \ti\om \, , \quad -(-\ti\om) = \ti\om \, , \quad \ga(-\ti\om) =
\ga(\ti\om) \quad \text{for } \ti\om \in \ti\Om \, ;
\]
these properties characterize the involution. In the case $ \ti\Om = \Om
\times \Om_\pm $ we have $ -(\om,s) = (\om,-s) $ for $ \om \in \Om $, $ s =
\pm1 $.

An isomorphism between two binary extensions boils down to an automorphism of
$ (\Om\times\Om_\pm,\ga) $. The general form of such automorphism is $ (\om,s)
\mapsto ( \om, s U(\om) ) $ for $ \om \in \Om $, $ s = \pm1 $; here $ U $ runs
over measurable functions $ \Om \to \{-1,+1\} $. The automorphism commutes
with the involution, thus, every isomorphism of extensions intertwines the
involutions,
\[
\theta (-\ti\om) = - \theta (\ti\om) \quad \text{for } \ti\om \in \ti\Om \, .
\]

\begin{definition}
(a) An \emph{inductive system of probability spaces} consists of probability
spaces $ \Om_n $ and measure preserving maps $ \be_n : \Om_n \to \Om_{n+1} $
for $ n = 1,2,\dots $

(b) Let $ (\Om_n,\be_n)_n $ and $ (\Om'_n,\be'_n)_n $ be two inductive systems
of probability spaces. A \emph{morphism} from $ (\Om_n,\be_n)_n $ to $
(\Om'_n,\be'_n)_n $ is a sequence of measure preserving maps $ \ga_n : \Om_n
\to \Om'_n $ such that the infinite diagram
\[
\xymatrix@C+1cm{
 \Om_1 \ar[d]^{\ga_1} \ar[r]^{\be_1} & \Om_2 \ar[d]^{\ga_2} \ar[r]^{\be_2} &
  \dots
\\
 \Om'_1 \ar[r]^{\be'_1} & \Om'_2 \ar[r]^{\be'_2} & \dots
}
\]
is commutative. If each $ \ga_n $ is invertible, the morphism is an
\emph{isomorphism.}

(c) A morphism $ (\ga_n)_n $ is \emph{binary,} if for every $ n $ the
extension $ (\Om_n,\ga_n) $ of $ \Om'_n $ is binary, and each $ \be_n $
intertwines the corresponding involutions,
\[
\be_n (-\om_n) = - \be_n (\om_n) \quad \text{for } \om_n \in \Om_n \, .
\]
Given a binary morphism $ (\ga_n)_n $ from $ (\Om_n,\be_n)_n $ to $
(\Om'_n,\be'_n)_n $, we say that $ (\Om_n,\be_n)_n $ is a \emph{binary
extension} of $ (\Om'_n,\be'_n)_n $ (according to $ (\ga_n)_n $).
\end{definition}

\begin{definition}
Let $ (\Om_n,\be_n)_n $ be an inductive system of probability spaces, $
(\ti\Om_n,\ti\be_n)_n $ its binary extension (according to $ (\ga_n)_n $), and
$ (\ti\Om'_n,\ti\be'_n)_n $ another binary extension of $ (\Om_n,\be_n)_n $
(according to $ (\ga'_n)_n $). An \emph{isomorphism} between the two binary
extensions is an isomorphism $ (\theta_n)_n $ between $ (\ti\Om_n,\ti\be_n)_n $
and $ (\ti\Om'_n,\ti\be'_n)_n $ treated as inductive systems of probability
spaces, satisfying the following condition: for each $ n $ the diagram
\[
\xymatrix{
 \ti\Om_n \ar[dr]_{\ga_n} \ar[rr]^{\theta_n} & & \ti\Om'_n \ar[dl]^{\ga'_n}
\\
 & \Om_n
}
\]
is commutative.
\end{definition}

In other words, an isomorphism between the two binary extensions of the
inductive system is a sequence $ (\theta_n)_n $ where each $ \theta_n $ is an
isomorphism between the two binary extensions $ (\ti\Om_n,\ti\ga_n) $ and $
(\ti\Om'_n,\ti\ga'_n) $ of the probability space $ \Om_n $, such that the
diagram
\[
\xymatrix@C+1cm{
 \ti\Om_n \ar[d]^{\theta_n} \ar[r]^{\ti\be_n} & \ti\Om_{n+1}
 \ar[d]^{\theta_{n+1}}
\\
 \ti\Om'_n \ar[r]^{\ti\be'_n} & \ti\Om'_{n+1}
}
\]
is commutative for every $ n $.

\begin{lemma}\label{5.43}
Let $ (\Om_n,\be_n)_n $ be an inductive system of probability spaces.

(a) Let $ X_n : \Om_n \to \{-1,+1\} $ be measurable functions, and
\begin{equation}\label{5.44}
\begin{gathered}
\ti\Om_n = \Om_n \times \Om_\pm \, , \\
\begin{gathered}
 \ti\be_n (\om_n,s_n) = \( \be_n(\om_n), s_n X_n(\om_n) \) \\
 \ga_n (\om_n,s_n) = \om_n
\end{gathered}
\quad \bigg\} \text{ for } \om_n \in \Om_n, \; s_n = \pm1 \, .
\end{gathered}
\end{equation}
Then $ (\ti\Om_n,\ti\be_n) $ is a binary extension of $ (\Om_n,\be_n)_n $
(according to $ (\ga_n)_n $).

(b) Every binary extension of $ (\Om_n,\be_n)_n $ is isomorphic to the
extension of the form \eqref{5.44}, for some $ (X_n)_n $.
\end{lemma}

\begin{proof}
(a) Clearly, $ \ti\be_n $ and $ \ga_n $ are measure preserving, $ \ga_n $ is
binary, and $ \ga_{n+1} (\ti\be_n(\om_n,s_n)) = \be_n(\om_n) = \be_n
(\ga_n(\om_n,s_n)) $.

(b) Let $ (\ti\Om_n,\ti\be_n)_n $ be a binary extension of $ (\Om_n,\be_n)_n $
according to $ (\ga_n)_n $. Without loss of generality we assume that $
\ti\Om_n = \Om_n \times \Om_\pm $ and $ \ga_n(\om_n,s_n) = \om_n $. The
relations $ \ga_{n+1} (\ti\be_n(\om_n,s_n)) = \be_n (\ga_n(\om_n,s_n)) =
\be_n(\om_n) $ and $ \ti\be_n(-\om_n) = -\ti\be_n(\om_n) $ show that $
\ti\be_n $ is of the form $  \ti\be_n (\om_n,s_n) = \( \be_n(\om_n), s_n
X_n(\om_n) \) $ for some measurable $ X_n : \Om_n \to \{-1,+1\} $.
\end{proof}

Given an inductive system $ (\Om_n,\be_n)_n $ of probability spaces and two
sequences $ (X_n)_n $, $ (Y_n)_n $ of measurable functions $ X_n, Y_n : \Om_n
\to \{-1,+1\} $, the construction \eqref{5.44} gives us two binary extensions
of $ (\Om_n,\be_n)_n $. One extension, $ (\ti\Om_n,\ti\be_n)_n $, $ (\ga_n)_n
$, corresponds to $ (X_n)_n $, the other extension, $ (\ti\Om'_n,\ti\be'_n)_n
$, $ (\ga'_n)_n $, corresponds to $ (Y_n)_n $. We want to know, whether they
are isomorphic or not.

For each $ n $ separately, the two binary extensions of the probability space
$ \Om_n $ coincide: $ \ti\Om_n = \Om_n \times \Om_\pm = \ti\Om'_n $, $
\ga_n(\om_n,s_n) = \om_n = \ga'_n(\om_n,s_n) $. Every isomorphism $ \theta_n $
between them is of the form
\[
\theta_n (\om_n,s_n) = \( \om_n, s_n U_n(\om_n) \) \quad \text{for } \om_n \in
\Om_n, \; s_n = \pm1 \, ,
\]
where $ U_n : \Om_n \to \{-1,+1\} $ is a measurable function. In order to form
an isomorphism between the binary extensions of the inductive system, these $
\theta_n $ must satisfy the condition $ \theta_{n+1} ( \ti\be_n(\ti\om_n) ) =
\ti\be'_n ( \theta_n(\ti\om_n) ) $, that is (recall \eqref{5.44}),
\[
X_n(\om_n) U_{n+1}(\be_n(\om_n)) = U_n(\om_n) Y_n(\om_n) \quad \text{for }
\om_n \in \Om_n \, .
\]

Given an inductive system $ (\Om_n,\be_n)_n $ of probability spaces, we
consider the commutative group $ G((\Om_n,\be_n)_n) $ of all sequences $ f =
(f_n)_n $ of measurable functions $ f_n : \Om_n \to \{-1,+1\} $ treated $
\modO $; the group operation is the pointwise multiplication. We define the
shift homomorphism $ T : G((\Om_n,\be_n)_n) \to G((\Om_n,\be_n)_n) $ by
\[
(Tf)_n (\om_n) = f_{n+1} \( \be_n(\om_n) \) \quad \text{for } \om_n \in \Om_n
\, .
\]
According to \eqref{5.44}, every $ X \in G((\Om_n,\be_n)_n) $ leads to a
binary extension of $ (\Om_n,\be_n)_n $. We summarize the previous paragraph
as follows.

\begin{lemma}\label{5.6}
The binary extensions corresponding to $ X,Y \in G((\Om_n,\be_n)_n) $ are
isomorphic if and only if $ X T(U) = Y U $ for some $ U \in G((\Om_n,\be_n)_n)
$.
\end{lemma}

\section[Binary extensions: Hilbert spaces]
 {\raggedright Binary extensions: Hilbert spaces}
\label{sec:4-5}
Given an extension of a probability space, $ \ga : \ti\Om \to \Om $, we have a
natural embedding of Hilbert spaces, $ L_2(\Om) \subset L_2(\ti\Om) $, and a
natural action of the commutative algebra $ L_\infty (\Om) $ on $ L_2(\ti\Om)
$. ($ L_2 $ and $ L_\infty $ over $ \C $ are meant.) Assume that the extension
is binary. Then the embedded subspace and its orthogonal complement are the
`even' and `odd' subspaces w.r.t.\ the involution $ \ti\om \mapsto -\ti\om $
on $ \ti\Om $; that is,
\begin{gather*}
\psi \in L_2(\Om) \quad \text{if and only if} \quad \psi(-\ti\om) =
 \psi(\ti\om) \text{ for almost all } \ti\om \in \ti\Om \, ; \\
\psi \in L_2(\ti\Om) \ominus L_2(\Om) \quad \text{if and only if} \quad
\psi(-\ti\om) = -\psi(\ti\om) \text{ for almost all } \ti\om \in \ti\Om \, .
\end{gather*}

\begin{lemma}\label{6.0}
Let $ \ga : \ti\Om \to \Om $ and $ \ga' : \ti\Om' \to \Om $ be two binary
extensions of a probability space $ \Om $. Then the following two conditions
on a unitary operator $ \Theta : L_2(\ti\Om') \to L_2(\ti\Om) $ are
equivalent:

(a) $ \Theta $ is trivial on $ L_2(\Om) $, and intertwines the two actions of
$ L_\infty(\Om) $. In other words,
\begin{gather*}
\Theta \psi = \psi \quad \text{for all } \psi \in L_2(\Om) \, , \\
\Theta ( h \cdot \psi ) = h \cdot ( \Theta \psi ) \quad \text{for all } \psi
 \in L_2(\ti\Om'), \; h \in L_\infty(\Om) \, .
\end{gather*}

(b) There exists an isomorphism of extensions $ \theta : \ti\Om \to \ti\Om' $
and $ h \in L_\infty(\Om) $, $ |h(\cdot)|=1 $, such that
\begin{gather*}
\Theta \psi = \psi \circ \theta \quad \text{for all } \psi \in L_2(\Om) \, ,
 \\
\Theta \psi = h \cdot ( \psi \circ \theta ) \quad \text{for all } \psi \in
 L_2(\ti\Om') \ominus L_2(\Om) \, .
\end{gather*}
\end{lemma}

\begin{proof}
(b) \imp (a): evident.
(a) \imp (b): Without loss of generality we assume that $ \ti\Om = \ti\Om' =
\Om \times \Om_\pm $ and $ \ga (\om,s) = \ga' (\om,s) = \om $. The Hilbert
space $ L_2(\ti\Om) \ominus L_2(\Om) $ consists of functions of the form $
(\om,s) \mapsto s f(\om) $ where $ f $ runs over $ L_2(\Om) $. Thus, $
L_2(\ti\Om) \ominus L_2(\Om) $ is naturally isomorphic to $ L_2(\Om) $, and
the isomorphism intertwines the actions of $ L_\infty(\Om) $. The operator $
\Theta $ maps $ L_2(\ti\Om') \ominus L_2(\Om) $ onto $ L_2(\ti\Om) \ominus
L_2(\Om) $ and leads to an operator $ L_2(\Om) \to L_2(\Om) $ that commutes
with $ L_\infty(\Om) $ and is therefore the multiplication by a function $ h
\in L_\infty(\Om) $.
\end{proof}

An inductive system of probability spaces $ (\Om_n,\be_n)_n $ leads evidently
to a decreasing sequence of Hilbert spaces,
\[
\xymatrix@C+1cm{
 L_2(\Om_1) & L_2(\Om_2) \ar@{_{(}->}[l] & \dots
 \ar@{_{(}->}[l]
}
\]
Similarly, a morphism from $ (\Om_n,\be_n)_n $ to $ (\Om'_n,\be'_n)_n $ leads
to a commutative diagram of Hilbert space embeddings
\[
\xymatrix@C+1cm{
 L_2(\Om_1) & L_2(\Om_2) \ar@{_{(}->}[l] & \dots
 \ar@{_{(}->}[l]
\\
 L_2(\Om'_1) \ar@{^{(}->}[u] & L_2(\Om'_2) \ar@{^{(}->}[u]
 \ar@{_{(}->}[l] & \dots \ar@{_{(}->}[l]
}
\]
The commutative algebra $ L_\infty (\Om'_n) $ acts on $ L_2 (\Om'_n) $ and $
L_2 (\Om_n) $, and the embedding $ L_2 (\Om'_n) \to L_2 (\Om_n) $ intertwines
these two actions.

\begin{lemma}\label{6.22}
Let $ (\Om_n,\be_n)_n $ be an inductive system of probability spaces, $
(\ti\Om_n,\ti\be_n)_n $ its binary extension (according to $ (\ga_n)_n $), and
$ (\ti\Om'_n,\ti\be'_n)_n $ another binary extension of $ (\Om_n,\be_n)_n $
(according to $ (\ga'_n)_n $). Then the following two conditions are
equivalent:

(a) The two binary extensions are isomorphic.

(b) There exist unitary operators
\[
\Theta_n : L_2 (\ti\Om'_n) \to L_2 (\ti\Om_n)
\]
such that for every $ n $, $ \Theta_n $ intertwines the actions of $
L_\infty(\Om_n) $ on $ L_2(\ti\Om_n) $ and $ L_2(\ti\Om'_n) $, and the
following two diagrams are commutative:
\[
\xymatrix{
 L_2(\ti\Om_n) & & L_2(\ti\Om'_n) \ar[ll]_{\Theta_n}
\\
 & L_2(\Om_n) \ar@{_{(}->}[ul] \ar@{^{(}->}[ur]
}
\qquad
\xymatrix{
 L_2(\ti\Om_n) & L_2(\ti\Om_{n+1}) \ar@{_{(}->}[l]
\\
 L_2(\ti\Om'_n) \ar@{^{(}->}[u]_{\Theta_n} & L_2(\ti\Om'_{n+1})
 \ar@{^{(}->}[u]_{\Theta_{n+1}} \ar@{_{(}->}[l]
}
\]
\end{lemma}

\begin{proof}
(a) \imp (b): evident.
(b) \imp (a):
For each $ n $ separately we have two binary extensions $ (\ti\Om_n,\ga_n) $,
$ (\ti\Om'_n,\ga'_n) $ of the probability space $ \Om_n $, and a unitary
operator $ \Theta_n : L_2(\ti\Om'_n) \to L_2(\ti\Om_n) $ that satisfies
Condition (a) of Lemma \ref{6.0}. On the other hand, due to Lemma \ref{5.43}
we may assume that $ \ti\Om_n = \Om_n \times \Om_\pm = \ti\Om'_n $, $
\ga_n(\om_n,s_n) = \om_n = \ga'_n(\om_n,s_n) $, $ \ti\be_n(\om_n,s_n) = \(
\be_n(\om_n), s_n X_n(\om_n) \) $ and $ \ti\be'_n(\om_n,s_n) = \(
\be_n(\om_n), s_n Y_n(\om_n) \) $.

Now Lemma \ref{6.0} gives us $ h_n \in L_\infty(\Om_n) $, $ |h_n(\cdot)|
= 1 $, such that $ \Theta_n \psi = h_n \cdot \psi $ for all $ \psi \in L_2 (
\Om_n \times \Om_\pm ) \ominus L_2(\Om_n) $. In other words, if $
\psi(\om_n,s_n) = s_n f(\om_n) $ then $ (\Theta_n \psi) (\om_n,s_n) = s_n
f(\om_n) h_n(\om_n) $; here $ f $ runs over $ L_2(\Om_n) $. By commutativity
of the second diagram, $ (\Theta_{n+1} \psi) \circ \ti\be_n = \Theta_n ( \psi
\circ \ti\be'_n ) $ for $ \psi \in L_2(\ti\Om'_{n+1}) $. For the case $
\psi(\om_{n+1},s_{n+1}) = s_{n+1} f(\om_{n+1}) $ we have, first,
\begin{multline*}
\( (\Theta_{n+1} \psi) \circ \ti\be_n \) (\om_n,s_n) = (\Theta_{n+1} \psi) \(
 \be_n(\om_n), s_n X_n(\om_n) \) = \\
= s_n X_n(\om_n) f(\be_n(\om_n)) h_{n+1} (\be_n(\om_n)) \, ,
\end{multline*}
and second,
\begin{gather*}
(\psi \circ \ti\be'_n) (\om_n,s_n) = \psi \( \be_n(\om_n), s_n Y_n(\om_n) \) =
 s_n Y_n(\om_n) f(\be_n(\om_n)) \, , \\
\Theta_n (\psi \circ \ti\be'_n) (\om_n,s_n) = s_n Y_n(\om_n) f(\be_n(\om_n))
h_n(\om_n) \, .
\end{gather*}
They are equal, which means that $ X_n(\om_n) h_{n+1} (\be_n(\om_n)) =
Y_n(\om_n) h_n(\om_n) $, that is,
\[
(h_{n+1} \circ \be_n) \cdot X_n = h_n \cdot Y_n \, .
\]
By Lemma \ref{5.6} it is sufficient to find measurable functions $ U_n : \Om_n
\to \{-1,+1\} $ such that
\[
(U_{n+1} \circ \be_n) \cdot X_n = U_n \cdot Y_n \quad \text{for all } n \, .
\]
We choose a Borel function $ \phi : \T \to \{-1,+1\} $, where $ \T = \{ z \in
\C : |z|=1 \} $, such that $ \phi(-z) = -\phi(z) $ for all $ z \in \T $. For
example, $ \phi(\E^{\I\al}) = +1 $ for $ \al \in [0,\pi) $ but $ -1 $ for $
\al \in [\pi,2\pi) $. The functions $ U_n(\cdot) = \phi(h_n(\cdot)) $ satisfy
the needed equation, since $ X_n(\cdot) = \pm1 $, $ Y_n(\cdot) = \pm1 $.
\end{proof}

\section[Products of binary extensions]
 {\raggedright Products of binary extensions}
\label{sec:4-8}
Definitions and statements of this section are used only in Sect.~\ref{sec:8}
(in the proof of Lemma \ref{11.7}).

Special measures are taken in the next definition in order to keep the product
binary (rather than quaternary).

\begin{definition}\label{4-8.1}
Let $ (\ti\Om_k,\ga_k) $ be a binary extension of a probability space $ \Om_k
$ for $ k=1,2 $; $ \Om = \Om_1 \times \Om_2 $; and $ A \subset \Om $ a
measurable set. The \emph{product} of these two binary extensions (according
to $ A $) is the extension $ (\ti\Om,\ga) $ of $ \Om $ defined as follows:
\[
\ti\Om = \underbrace{ \{ ( \ti\om_1, \om_2 ) : ( \ga_1(\ti\om_1),\om_2 ) \in A
\} }_{ \ti A } \uplus \underbrace{ \{ (\om_1,\ti\om_2) : (
\om_1,\ga_2(\ti\om_2) ) \in \Om \setminus A \} }_{ \ti\Om \setminus \ti A } \,
,
\]
the measure on $ \ti A $ is induced from (the product measure on) $ \ti\Om_1
\times \Om_2 $, on $ \ti\Om \setminus \ti A $ --- from $ \Om_1 \times \ti\Om_2
$;
\[
\ga(\ti\om_1,\om_2) = (\ga_1(\ti\om_1),\om_2) \, , \quad
\ga(\om_1,\ti\om_2) = (\om_1,\ga_2(\ti\om_2)) \, .
\]
\end{definition}

Here and henceforth $ \om_k $ runs over $ \Om_k $, and $ \ti\om_k $ runs over
$ \ti\Om_k $.

Clearly, the extension $ (\ti\Om,\ga) $ is binary.

Let a binary extension $ (\ti\Om,\ga) $ of $ \Om = \Om_1 \times \Om_2 $ be the
product of two binary extensions $ (\ti\Om_k,\ga_k) $, $ k=1,2 $ (according to
a given $ A \subset \Om $). Then we have a natural embedding of Hilbert
spaces,
\begin{equation}\label{4-8.*}
L_2(\ti\Om) \subset L_2(\ti\Om_1) \otimes L_2(\ti\Om_2) \, ;
\end{equation}
it arises from the natural measure preserving map $ \ti\Om_1 \times \ti\Om_2
\to \ti\Om $,
\[
(\ti\om_1,\ti\om_2) \mapsto \begin{cases}
 (\ti\om_1,\ga_2(\ti\om_2)) &\text{if } (\ga_1(\ti\om_1),\ga_2(\ti\om_2)) \in
  A, \\
 (\ga_1(\ti\om_1),\ti\om_2) &\text{otherwise}.
\end{cases}
\]
The restriction of the embedding \eqref{4-8.*} to $ L_2(\Om) $ is just the
tensor product of the two embeddings $ L_2(\Om_k) \subset L_2(\ti\Om_k) $, $
k=1,2 $, since the corresponding composition map $ \ti\Om_1 \times \ti\Om_2
\to \ti\Om \to \Om $ is just $ \ga_1 \times \ga_2 $.

The projection map $ \ti A \to \ti\Om_1 $, $ (\ti\om_1,\om_2) \mapsto \ti\om_1
$, need not be measure preserving, but anyway, generates a sub-\sif\ $ \F_1 $
on $ \ti A $.

\begin{lemma}\label{4-8.2}
Let $ (\ti\Om_k,\ga_k) $ and $ (\ti\Om'_k,\ga'_k) $ be two binary extensions
of a probability space $ \Om_k $ (for $ k=1,2 $), $ \Om = \Om_1 \times \Om_2 $,
$ A \subset \Om $ a measurable set, $ \Theta_k : L_2(\ti\Om_k) \to
L_2(\ti\Om'_k) $ unitary operators, each satisfying Condition (a) of Lemma
\ref{6.0}. Then $ \Theta_1 \times \Theta_2 $ maps $ L_2(\ti\Om) $ onto $
L_2(\ti\Om') $, $ L_2(\ti A) $ onto $ L_2(\ti A') $, and $ L_2(\ti A,\F_1) $
onto $ L_2(\ti A', \F'_1) $.
\end{lemma}

It is meant that $ L_2(\ti A,\F_1) \subset L_2(\ti A) \subset L_2(\ti A)
\oplus L_2 (\ti\Om \setminus \ti A) = L_2(\ti\Om) \subset L_2(\ti\Om_1)
\otimes L_2(\ti\Om_2) $ and $ L_2(\ti A',\F'_1) \subset L_2(\ti A') 
\oplus L_2 (\ti\Om' \setminus \ti A') = L_2(\ti\Om') \subset L_2(\ti\Om'_1)
\otimes L_2(\ti\Om'_2) $.

The reader may prove Lemma \ref{4-8.2} via Lemma \ref{6.0}, but the proof
below does not use Lemma \ref{6.0}.

\begin{proof}
The operator $ \Theta = \Theta_1 \otimes \Theta_2 $ intertwines the actions of
$ L_\infty (\Om_1) $ and $ L_\infty (\Om_2) $, therefore, also the actions of
$ L_\infty (\Om_1 \times \Om_2) $. In particular,
\[
\Theta \One_{\ti A} \psi = \One_{\ti A'} \Theta \psi \quad \text{for } \psi
\in L_2 (\ti\Om_1 \times \ti\Om_2) \, .
\]
The space $ L_2(\ti A) $ is the closure of linear combinations of vectors of
the form $ \psi = \One_{\ti A} (\phi_1 \otimes \phi_2) $, where $ \phi_1 \in
L_2(\ti\Om_1) $ and $ \phi_2 \in L_2(\Om_2) $. For such $ \psi $ we have
\[
\Theta \psi = \Theta \One_{\ti A} (\phi_1 \otimes \phi_2) = \One_{\ti A'}
\Theta (\phi_1 \otimes \phi_2) = \One_{\ti A'} (\Theta_1 \phi_1 \otimes
\Theta_2 \phi_2) \in L_2(\ti A') \, ,
\]
since $ \Theta_1 \phi_1 \in L_2(\ti\Om'_1) $ and $ \Theta_2 \phi_2 = \phi_2
\in L_2(\Om_2) $. Therefore $ \Theta(L_2(\ti A)) \subset L_2(\ti A') $. The
special case $ \phi_2 = \One $ gives $ \Theta(L_2(\ti A,\F_1)) \subset L_2(\ti
A',\F'_1) $. The same holds for $ \Theta^{-1} $, thus, the inclusions are in
fact equalities. Similarly, $ \Theta(L_2(\ti\Om \setminus \ti A)) =
L_2(\ti\Om' \setminus\ti A') $. It follows that $ \Theta(L_2(\ti\Om)) =
L_2(\ti\Om') $.
\end{proof}

\section[Some necessary conditions of isomorphism]
 {\raggedright Some necessary conditions of isomorphism}
\label{sec:5}
Let $ \mu_1 $ be a probability measure on the space $ \R^\infty $ (of all
infinite sequences of reals), $ \be : \R^\infty \to \R^\infty $ the shift, $
\be(x_1,x_2,\dots) = (x_2,x_3,\dots) $, and $ \mu_n $ the image of $ \mu_1 $
under $ \be^{n-1} $. Probability spaces $ \Om_n = (\R^\infty,\mu_n) $ with
maps $ \be_n = \be $ are an inductive system of probability spaces.

Let Borel functions $ f_n : \R \to \{-1,+1\} $ be given. We define $ X_n :
\Om_n \to \{-1,+1\} $ by
\[
X_n (x_n,x_{n+1},\dots) = f_n (x_n)
\]
and consider the corresponding binary extension of $ (\Om_n,\be_n)_n
$. Another sequence of functions $ g_n : \R \to \{-1,+1\} $ leads to another
binary extension. According to Lemma \ref{5.6} the two binary extensions are
isomorphic if and only if there exist $ U_n : \Om_n \to \{-1,+1\} $ such that
\begin{equation}\label{5.1}
U_{n+1} (x_{n+1},x_{n+2},\dots) = U_n (x_n,x_{n+1},\dots) f_n(x_n) g_n(x_n) \, .
\end{equation}
Functions that do not depend on $ x_n $, that is, functions of the form
\[
(x_1,x_2,\dots) \mapsto \phi ( x_1,\dots,x_{n-1}, x_{n+1},x_{n+2},\dots)
\]
are a subspace $ H_n \subset L_2(\mu_1) $. We consider vectors $ \psi_n \in
L_2(\mu_1) $,
\begin{equation}\label{psi}
\psi_n (x_1,x_2,\dots) = f_n (x_n) g_n (x_n) \, ,
\end{equation}
and the distance between $ \psi_n $ and $ H_n $.

\begin{lemma}\label{5.2}
The condition
\[
\dist ( \psi_n, H_n ) \to 0 \quad \text{as } n \to \infty
\]
is necessary for the two binary extensions to be isomorphic.
\end{lemma}

\begin{proof}
Let $ U_n $ satisfy \eqref{5.1}, then
\[
U_n (x_n,x_{n+1},\dots) = U_1 (x_1,x_2,\dots) h_1(x_1) \dots h_{n-1}(x_{n-1})
\, ,
\]
where $ h_n(x) = f_n(x) g_n(x) $. We have
\begin{multline*}
\psi_n (x_1,x_2,\dots) = h_n(x_n) = U_n (x_n,x_{n+1},\dots) U_{n+1}
 (x_{n+1},x_{n+2},\dots) = \\
= U_1 (x_1,x_2,\dots) h_1(x_1) \dots h_{n-1}(x_{n-1}) U_{n+1}
 (x_{n+1},x_{n+2},\dots) \, .
\end{multline*}
Taking into account that $ H_n $ is invariant under multiplication by any
(bounded measurable) function of $ x_1,\dots,x_{n-1} $ and $ x_{n+1}, x_{n+2},
\dots $ we see that $ \dist ( \psi_n, H_n ) \linebreak[0]
= \dist ( U_1, H_n ) $. The latter
converges to $ 0 $, since $ H_n $ contains all functions of $
x_1,\dots,x_{n-1} $.
\end{proof}

The conditional distribution of $ x_n $ given $ x_1,\dots,x_{n-1} $ and $
x_{n+1}, x_{n+2}, \dots $ (assuming that $ (x_1,x_2,\dots) $ is distributed $
\mu_1 $) is a probability measure $ \nu_n $ on $ \R $; this $ \nu_n $ is
random in the sense that it depends on $ x_1,\dots,x_{n-1} $ and $ x_{n+1},
x_{n+2}, \dots $ (whose distribution is the marginal of $ \mu_1 $).

Here is a useful condition on $ \mu_1 $:
\begin{equation}\label{5.4}
\exists \eps > 0 \;\; \forall n \;\; \Pr{ \text{$ \nu_n $ is $ \eps $-good} }
\ge \eps \, ,
\end{equation}
where a probability measure $ \nu $ on $ \R $ is called $ \eps $-good if
\begin{equation}\label{5.45}
\exists x \in \R \;\; \forall A \;\;\; \nu(A) \ge \eps \mes \( A \cap
(x,x+\eps) \)
\end{equation}
($ A $ runs over Borel subsets of $ \R $; `$ \mes $' stands for Lebesgue
measure).

Usually $ \nu $ has a density; then \eqref{5.45} requires the density to
exceed $ \eps $ on some interval of length $ \eps $.

\begin{lemma}\label{5.5}
Let $ \mu_1 $ satisfy \eqref{5.4}, and numbers $ \eps_n \in (0,\infty) $
satisfy $ \eps_n \to 0 $. Then there exist Borel functions $ f_n : \R \to
\{-1,+1\} $ such that for every $ c \in \R \setminus \{0\} $, defining $ g_n :
\R \to \{-1,+1\} $ by $ g_n(x) = f_n (x+c\eps_n) $ we get $ \psi_n $ (see
\eqref{psi}) violating the necessary condition of Lemma \ref{5.2} (and
therefore, $ (f_n)_n $ and $ (g_n)_n $ lead to two nonisomorphic binary
extensions).
\end{lemma}

\begin{proof}
We take $ \la_1, \la_2, \dots $ such that
\[
\la_n \eps_n = \begin{cases}
 1 &\text{for $ n $ odd}, \\
 \sqrt2 &\text{for $ n $ even},
\end{cases}
\]
and define
\[
f_n (x) = \si ( \la_n x ) \, ,
\]
where
\[
\si(x) = \begin{cases}
 -1 &\text{for } x \in \cup_{k\in\Z} [k-0.5,k) \, , \\
 +1 &\text{for } x \in \cup_{k\in\Z} [k,k+0.5) \, .
\end{cases}
\]
Let $ c \in \R $ be given, $ c \ne 0 $. It is sufficient to prove that at
least one of two claims
\[
\limsup_n \dist ( \psi_{2n}, H_{2n} ) > 0 \, , \quad \limsup_n \dist (
\psi_{2n-1}, H_{2n-1} ) > 0
\]
holds. Here $ \psi_n (x_1,x_2,\dots) = h_n(x_n) = f_n(x_n) g_n(x_n) = f_n(x_n)
f_n(x_n+c\eps_n) = \si (\la_n x_n) \si (\la_n x_n + c \la_n \eps_n ) $. The
function $ h_n $ is periodic, with period $ 1/\la_n $. The mean value $ M_n $
of $ h_n $ over the period is
\[
M_n = \la_n \int_0^{1/\la_n} h_n(x) \, \D x = \int_0^1 \si(u)
\si(u+c\la_n\eps_n) \, \D u \, .
\]
It reaches $ \pm1 $ when $ 2c\la_n\eps_n \in \Z $; otherwise $ -1 < M_n < 1
$. The relations $ 2c \in \Z $ and $ 2c \sqrt2 \in \Z $ are incompatible,
therefore at least one of two claims
\[
\sup_n |M_{2n}| < 1 \, , \quad \sup_n |M_{2n-1}| < 1
\]
holds. (Of course, $ M_{2n} $ and $ M_{2n-1} $ do not depend on $ n $, but
this fact does not matter.) It is sufficient to prove that
\begin{align*}
\sup_n |M_{2n}| < 1 \quad &\text{implies} \quad \limsup_n \dist ( \psi_{2n},
 H_{2n} ) > 0 \, , \\
\sup_n |M_{2n-1}| < 1 \quad &\text{implies} \quad \limsup_n \dist (
\psi_{2n-1}, H_{2n-1} ) > 0 \, .
\end{align*}
The former implication will be proved (the latter is similar). Assume the
contrary: $ \sup_n |M_{2n}| < 1 $ and $ \dist ( \psi_{2n}, H_{2n} ) \to 0 $.

For any probability measure $ \nu $ on $ \R $, the squared distance in the
space $ L_2(\nu) $ between the function $ h_n $ and the one-dimensional space
of constant functions is
\[
\int \bigg( h_n - \int h_n \, \D\nu \bigg)^2 \, \D\nu = \int h_n^2 \, \D\nu -
\bigg( \int h_n \, \D\nu \bigg)^2 = 1 - \bigg( \int h_n \, \D\nu \bigg)^2 \, .
\]
We use the random measure $ \nu_n $, take the average and recall the
definition of $ H_n $:
\[
\Ex \bigg( 1 - \bigg( \int h_n \, \D\nu_n \bigg)^2 \bigg) = \dist^2
(\psi_n,H_n) \, .
\]
Taking into account that $ \dist ( \psi_{2n}, H_{2n} ) \to 0 $ we see that $ |
\int h_{2n} \, \D\nu_{2n} | \to 1 $ in probability. In order to get a
contradiction to \eqref{5.4} it is sufficient to prove that $ \limsup_n
\sup_\nu | \int h_{2n} \, \D\nu | < 1 $, where $ \nu $ runs over all $ \eps
$-good measures (recall \eqref{5.4} and \eqref{5.45}). Or, equivalently,
\[
\liminf_n \inf_\nu \nu \( h_{2n}^{-1} (-1) \) > 0 \, , \quad \liminf_n \inf_\nu \nu \( h_{2n}^{-1} (+1) \) > 0 \, .
\]
The former will be proved (the latter is similar). By \eqref{5.45}, $ \nu \(
h_{2n}^{-1} (-1) \) \ge \eps \mes \( h_{2n}^{-1} (-1) \cap (x,x+\eps) \)
$. For large $ n $ the period $ 1/\la_{2n} $ of the function $ h_{2n} $ is $
\ll\eps $, therefore
\begin{multline*}
\frac1\eps \mes \( h_{2n}^{-1} (-1) \cap (x,x+\eps) \) \ge \frac12 \cdot
 \la_{2n} \int_0^{1/\la_{2n}} \frac{1-h_{2n}(x)}{2} \, \D x = \\
= \frac12 \cdot \frac{ 1 - M_{2n} }{ 2 } \ge \frac14 \( 1 - \sup_n |M_{2n}| \)
> 0 \, .
\end{multline*}
\end{proof}

\begin{remark}
The functions $ f_n $ constructed in the proof of Lemma \ref{5.5} depend only
on the numbers $ \eps_n $, not on the measure $ \mu_1 $.
\end{remark}

Let a probability measure $ \mu $ on $ C[0,1] $ be given. Random variables $
A(t) $ on the probability space $ \Om = (C[0,1],\mu) $ defined by $ A(t)(a) =
a(t) $ for $ a \in C[0,1] $, $ t \in [0,1] $, are a random process. For every
$ n $ the restriction map $ C[0,1] \to C[0,3^{-n}] $ sends $ \mu $ to some $
\mu_n $. Probability spaces $ \Om_n = (C[0,3^{-n}],\mu_n) $ with restriction
maps are an inductive system.

Given Borel functions $ f_n : \R \to \{-1,+1\} $, we define random variables $
X_n : \Om_n \to \{-1,+1\} $ by $ X_n = f_n ( A(2\cdot 3^{-n-1}) ) $. The
corresponding binary extension may be visualized as follows. We consider
pairs $ (a,s) $ of a function $ a \in C[0,1] $ and another function $ s :
(0,1] \to \{-1,+1\} $ constant on each $ [2\cdot 3^{-n-1},2\cdot 3^{-n}) $ and
such that $ s(2\cdot 3^{-n}-) s(2\cdot 3^{-n}) = f_{n-1} (a(2\cdot 3^{-n})) $
for all $ n $. We get a pair of random processes $ A(\cdot) $, $ S(\cdot) $
satisfying
\[
\frac{ S(2\cdot 3^{-n}) }{ S(2\cdot 3^{-n}-) } = f_{n-1} (A(2\cdot 3^{-n})) \,
.
\]
Their restrictions to $ [0,3^{-n}] $ give $ \ti\Om_n $. For each $ t $
(separately), the random variable $ S(t) $ is independent of the process $
A(\cdot) $ and takes on the two equiprobable values $ \pm1 $.

As before, given also $ g_n : \R \to \{-1,+1\} $ (thus, another binary
extension), we define $ \psi_n \in L_2(\Om) $ by
\[
\psi_n = f_n (A(2\cdot 3^{-n-1})) g_n (A(2\cdot 3^{-n-1})) \, .
\]
We consider the subspaces $ H_n \subset L_2(\Om) $ consisting of functions of
$ A(t) $ for $ t \in [0,3^{-n-1}] \cup [3^{-n},1] $ only (in other words,
functions of the restrictions of sample paths to $ [0,3^{-n-1}] \cup
[3^{-n},1] $).

\begin{lemma}\label{5.lenew}
The condition $ \dist(\psi_n,H_n) \to 0 $ is necessary for the two binary
extensions to be isomorphic.
\end{lemma}

The proof, similar to the proof of Lemma \ref{5.2}, is left to the reader.

Similarly, $ C[-1,1] $ may be used (instead of $ C[0,1] $), with $ \Om_n = (
C[-1,3^{-n}], \linebreak[0]
\mu_n) $; the process $ S(\cdot) $ jumps at $ 2 \cdot 3^{-n} $,
as before. Now $ H_n $ consists of functions of the restriction $ A
|_{[-1,3^{-n-1}]\cup[3^{-n},1]} $ of $ A $ to $ [-1,3^{-n-1}] \cup [3^{-n},1]
$ (rather than $ [0,3^{-n-1}] \cup [3^{-n},1] $). Lemma \ref{5.lenew} remains
true.

The conditional distribution of $ A(2\cdot3^{-n-1}) $ given $ A
|_{[-1,3^{-n-1}]\cup[3^{-n},1]} $ is too concentrated (when $ n $ is large)
for being $ \eps $-good (recall \eqref{5.45}). A useful condition on $ \mu $
stipulates rescaling by $ 3^{n/2} $:
\begin{gather}
\text{there exists $ \eps>0 $ such that for every $ n $,} \label{5.4new} \\
\text{the conditional distribution of $ 3^{n/2} A(2\cdot3^{-n-1}) $ given $ A
 |_{[-1,3^{-n-1}]\cup[3^{-n},1]} $} \notag \\
\text{is $ \eps $-good with probability $ \ge \eps $.} \notag
\end{gather}
Here is a counterpart of Lemma \ref{5.5} for $ \eps_n = 3^{-n/2} $.

\begin{lemma}\label{5.5new}
Let $ \mu $ satisfy \eqref{5.4new}. Then there exist Borel
functions $ f_n : \R \to \{-1,+1\} $ such that for every $ c \in \R \setminus
\{0\} $, defining $ g_n : \R \to \{-1,+1\} $ by $ g_n(x) = f_n (x+3^{-n} c) $
we get $ \psi_n = f_n (A(2\cdot 3^{-n-1})) g_n (A(2\cdot 3^{-n-1})) $
violating the necessary condition of Lemma \ref{5.lenew} (and
therefore, two nonisomorphic binary extensions).
\end{lemma}

The proof, similar to the proof of Lemma \ref{5.5}, is left to the reader.

\section[A binary extension of Brownian motion]
 {\raggedright A binary extension of Brownian motion}
\label{sec:6}
The space $ C[0,1] $ of all continuous functions $ b : [0,1] \to \R $, endowed
with the Wiener measure $ \W $, is a probability space. Random variables $
B(t) $ on $ (C[0,1],\W) $, defined for $ t \in [0,1] $ by $ B(t)(b) = b(t) $,
are the Brownian motion on $ [0,1] $. Almost surely, a Brownian sample path on
$ [0,1] $ has a unique (global) minimum,
\[
\min_{t\in[0,1]} B(t) = B(\tau) \, ,
\]
$ \tau $ being a measurable function on $ (C[0,1],\W) $, $ 0<\tau(\cdot)<1 $
a.s.

We define another random process $ A $, on the time interval $ [-1,1] $, by
\begin{align*}
A(t) &= B \( \min(1,\tau+t)\) - B(\tau) \quad \text{for } t \in [0,1] \, , \\
A(t) &= B \( \max(0,\tau+t)\) - B(\tau) \quad \text{for } t \in [-1,0] \, .
\end{align*}
A \measurable{\W} map $ C[0,1] \to C[-1,1] $ is thus introduced. The map is
one-to-one ($ \modO $), since $ B(\cdot) $ is non-constant on every time
interval, almost surely.

\begin{proposition}\label{prop8}
The process $ A $ satisfies \eqref{5.4new}.
\end{proposition}

The proof is given after three lemmas.

The conditional distribution of the process $ B $ given the restriction $
A|_{[-1,\eps]} $ (for a given $ \eps \in (0,1) $) is the same as the
conditional distribution of the process $ B $ given $ \tau $ and $
B|_{[0,\tau+\eps]} $, since the two corresponding measurable partitions of $
(C[0,1],\W) $ are equal ($ \modO $). This conditional distribution is a
probability measure on the set of Brownian sample paths $ b $ such that
\begin{align*}
b(t) = x(t) &\quad \text{for } t \in [0,s+\eps] \, , \\
b(t) > x(s) &\quad \text{for } t \in [s+\eps,1] \, ;
\end{align*}
here $ s \in (0,1-\eps) $ is a given value of $ \tau $, and $ x \in
C[0,s+\eps] $ is a given sample path of $ B|_{[0,\tau+\eps]} $; of course,
$ s $ is the unique minimizer of $ x $. We assume that $ s < 1-\eps $, since
the other case is trivial (the conditional distribution is a single atom).

The corresponding conditional distribution of $ B|_{[s+\eps,1]} $ is a
probability measure on the set of functions $ b \in C[s+\eps,1] $ such that
\begin{align*}
& b(s+\eps) = x(s+\eps) \, , \\
& b(t) > x(s) \quad \text{for } t \in [s+\eps,1] \, .
\end{align*}
This set depends only on the three numbers $ s+\eps $, $ x(s+\eps) $, and $
x(s) $. One may guess that the considered measure on this set also depends on
these three numbers only (rather than the whole function $ x $). The following
well-known lemma confirms the guess and gives a simple description of the
measure.

\begin{lemma}\label{6.1}
The conditional distribution of $ B|_{[s+\eps,1]} $ given that $ \tau = s $
and $ B_{[0,\tau+\eps]} = x $ is equal to the conditional distribution of $
B|_{[s+\eps,1]} $ given that $ B(s+\eps) = x(s+\eps) $ and $ B(t) > x(s) $ for
$ t \in [s+\eps,1] $.
\end{lemma}

\begin{proof}
We take $ n $ such that $ \frac1n < \eps $. Let $ k \in \{1,\dots,n-1\} $. The
conditional distribution of $ B|_{[\frac k n,1]} $ given $ B|_{[0,\frac k n]}
$ depends only on $ B\(\frac k n\) $ (by the Markov property of $ B $) and is
just the distribution of the Brownian motion starting from $ \( \frac k n, B\(
\frac k n \) \) $. Therefore the conditional distribution of $ B|_{[\frac k
n,1]} $ given both $ B|_{[0,\frac k n]} $ and $ \frac{k-1}n < \tau < \frac k n
$ is the distribution of the Brownian motion starting from $ \( \frac k n, B\(
\frac k n \) \) $ and conditioned to stay above the minimum of the given path
on $ [0,\frac k n] $. (Indeed, a measurable partition of the whole probability
space induces a measurable partition of a given subset of positive
probability, and conditional measures for the former partition induce
conditional measures for the latter partition.) Now it is easy to condition
further on $ B|_{[\frac k n,\tau+\eps]} $ and combine all $ k $ together.
\end{proof}

Lemma \ref{6.1} gives the conditional distribution of $ B|_{[\tau+\eps,1]} $
given $ \tau $ and $ B|_{[0,\tau+\eps]} $. Now we turn to the conditional
distribution of $ B|_{[\tau+\eps,\tau+3\eps]} $ given $ \tau $, $
B|_{[0,\tau+\eps]} $ and $ B|_{[\tau+3\eps,1]} $ (in the case $ \tau+3\eps<1
$). We are especially interested in $ B(\tau+2\eps) $.

\begin{lemma}\label{6.2}
The conditional distribution of $ B(\tau+2\eps) - B(\tau) $, given $ \tau $
(such that $ \tau+3\eps < 1 $), $ B|_{[0,\tau+\eps]} $ and $
B|_{[\tau+3\eps,1]} $, has the density
\[
x \mapsto \frac{ \( 1 - \E^{-2ax/\eps} \) \( 1 - \E^{-2bx/\eps} \) }{ 1 -
\E^{-ab/\eps} } \cdot \frac1{ \sqrt{\pi\eps} } \exp \bigg( \! -\frac1\eps \Big(
x - \frac{a+b}2 \Big)^2 \bigg) \quad \text{for } x > 0 \, ,
\]
where $ a = B(\tau+\eps) - B(\tau) $ and $ b = B(\tau+3\eps) - B(\tau) $.
\end{lemma}

\begin{proof}
Using Lemma \ref{6.1} we turn to an equivalent question: a Brownian motion
starting from $ \( s+\eps, B(s+\eps) \) $ is conditioned to stay above $ B(s)
$ on $ [s+\eps,s+3\eps] $, and is known on $ [s+3\eps,1] $ (which means
another conditioning, of course); we need the (conditional) distribution of $
B(s+2\eps) $. Omitting for a while the condition $ B_{[s+\eps,s+3\eps]}(\cdot)
> B(s) $ we get the so-called Brownian bridge, --- the Brownian motion on $
[s+\eps,s+3\eps] $ with given boundary values $ B(s+\eps) $, $ B(s+3\eps)
$. (Later we'll condition the bridge to stay above $ B(s) $.)

For the bridge, $ B(s+2\eps) $ has the normal distribution $ N \( \frac{
B(s+\eps) + B(s+3\eps) }2, \frac\eps2 \) $. Given $ B(s+2\eps) $ we get two
independent bridges, one on $ [s+\eps,s+2\eps] $, the other on $
[s+2\eps,s+3\eps] $. The bridge on $ [s+\eps,s+2\eps] $ stays above $ B(s) $
with the probability (calculated via the reflection principle)
\[
\frac{ p_\eps (a-x) - p_\eps (a+x) }{ p_\eps (a-x) } = 1 - \exp \Big( \! -
 \frac2\eps ax \Big) \, ,
\]
where $ a = B(s+\eps) - B(s) $, $ x = B(s+2\eps) - B(s) $, $ b = B(s+3\eps) -
B(s) $, and
\[
p_\eps(u) = \frac1{\sqrt{2\pi\eps}} \exp \Big( -\frac{u^2}{2\eps} \Big) \, .
\]
It remains to write similar formulas on $ [s+2\eps,s+3\eps] $ and the whole $
[s+\eps,s+3\eps] $, and apply the Bayes formula
\[
p_{X|A} (x) = \frac{ \cP{A}{X=x} p_X(x) }{ \Pr{A} }
\]
for the conditional density $ p_{X|A} (\cdot) $ of a random variable $ X $
given an event $ A $. Namely, $ X = B(s+2\eps) - B(s) \sim N \( \frac{a+b}2,
\frac\eps2 \) $, $ p_X(x) = p_{\eps/2} \( x - \frac{a+b}2 \) $, $ A $ is the
event $ B_{[s+\eps,s+3\eps]}(\cdot) > B(s) $, $ \cP{A}{X=x} = \( 1 -
\E^{-2ax/\eps} \) \( 1 - \E^{-2bx/\eps} \) $, $ \Pr A = 1 - \E^{-ab/\eps} $.
\end{proof}

\begin{lemma}\label{6.3}
There exists $ \eps > 0 $ such that for all $ a,b \in (0,\infty) $ the
probability measure on $ (0,\infty) $ that has the density
\[
x \mapsto \frac{ \( 1 - \E^{-2ax} \) \( 1 - \E^{-2bx} \) }{ 1 - \E^{-ab} }
\cdot \frac1{ \sqrt{\pi} } \exp \Big( - \Big( x - \frac{a+b}2 \Big)^2 \Big)
\]
is $ \eps $-good (as defined by \eqref{5.45}).
\end{lemma}

\begin{proof}
It is sufficient to prove that 
\[
\inf_{a,b\in (0,\infty)} \inf_{x \in [\frac{a+b}2+1,\frac{a+b}2+2]} p_{a,b}(x)
> 0 \, ,
\]
where $ p_{a,b}(\cdot) $ is the given density. Assume the contrary: there
exist $ a_n, b_n, x_n $ such that $ b_n \ge a_n > 0 $, $ \frac{a_n+b_n}2+1 \le
x_n \le \frac{a_n+b_n}2+2 $, and $ p_{a_n,b_n} (x_n) \to 0 $. Then
\[
\frac{ \( 1 - \E^{-2a_n x_n} \) \( 1 - \E^{-2b_n x_n} \) }{ 1 - \E^{-a_n b_n}
} \to 0 \, .
\]
It follows that $ 1 - \E^{-2a_n x_n} \to 0 $, $ a_n x_n \to 0 $, $ a_n \(
\frac{a_n+b_n}2+1 \) \to 0 $, $ a_n \to 0 $, $ a_n b_n \to 0 $, and
\[
\frac{ \( 1 - \E^{-2a_n x_n} \) \( 1 - \E^{-2b_n x_n} \) }{ a_n b_n } \to 0 \,
.
\]
Therefore
\[
\frac{ 1 - \E^{-2a_n} }{ a_n } \cdot \frac{ 1 - \E^{-2b_n} }{ b_n } \to 0 \, ;
\quad \frac{ 1 - \E^{-2b_n} }{ b_n } \to 0 \, ;
\]
$ b_n \to \infty $, $ b_n x_n \to \infty $; $ 1 - \E^{-2b_n x_n} \to 1 $;
\[
\frac{ 1 - \E^{-2a_n x_n} }{ a_n b_n } \to 0 \, ; \quad \frac{ 1 - \E^{-a_n
b_n} }{ a_n b_n } \to 0 \, ,
\]
in contradiction to $ a_n b_n \to 0 $.
\end{proof}

\begin{proof}[Proof of Prop.~\ref{prop8}]
Lemma \ref{6.2} (for $ \eps = 3^{-n-1} $) gives us the conditional
distribution of $ A ( 2 \cdot 3^{-n-1} ) $ given $ \tau $ and $
A|_{[-1,3^{-n-1}]\cup[3^{-n},1]} $, but only for the case $ \tau + 3^{-n} < 1
$. Lemma \ref{6.3} states that the corresponding distribution of $ 3^{n/2} A (
2 \cdot 3^{-n-1} ) $ is $ \eps $-good. It remains to note that $ \Pr{ \tau +
3^{-n} < 1 } \ge \Pr{ \tau < 2/3 } \ge \eps $.
\end{proof}

Combining Prop.~\ref{prop8} and Lemma \ref{5.5new} we get a binary extension
of the inductive system (of probability spaces) formed by the restrictions $
A|_{[-1,3^{-n}]} $ of the process $ A $. In terms of the Brownian motion $ B $
the inductive system is formed by $ B|_{[0,\tau+3^{-n}]} $, and the binary
extension may be visualized by a random function $ S : (\tau,1) \to \{-1,+1\}
$ constant on $ [\tau+2\cdot3^{-n-1}, \tau+2\cdot3^{-n}) \cap (0,1) $ for each
$ n $ and such that
\begin{equation}\label{6.9}
\frac{ S(\tau+2\cdot3^{-n}) }{ S(\tau+2\cdot3^{-n}-) } = f_{n-1} \(
B(\tau+2\cdot3^{-n}) - B(\tau) \)
\end{equation}
for all $ n $ such that $ \tau+2\cdot3^{-n} < 1 $. Here $ f_n : \R \to
\{-1,+1\} $ are the functions given by Lemma \ref{5.5new}. They are
constructed as to make the binary extension \emph{sensitive to drift} in the
following sense. For every $ c \in \R \setminus \{0\} $ the binary extension
constructed via
\[
f_n \( B(\tau+2\cdot3^{-n}) - B(\tau) + c \cdot 2\cdot3^{-n} \)
\]
is not isomorphic to that for $ c=0 $.

\section[A new noise extending the white noise]
 {\raggedright A new noise extending the white noise}
\label{sec:7}
\smallskip

\parbox{6.2cm}{\small\textit{%
This is a noise richer than white noise: in addition to the increments of a
Brownian motion $ B $ it carries a countable collection of independent
Bernoulli random variables which are attached to the local minima of $ B $.}\\
\mbox{}\hfill J.~Warren \cite[the end]{Wa1}
}\hfill
\parbox{6.9cm}{\small\textit{%
\dots magically, this independent random variable has appeared from
somewhere! Indeed, it really has appeared from thin air, because\dots it is
not present at time $ 0 $!}\\ 
\mbox{}\hfill L.C.G.~Rogers, D.~Williams \cite[p.~156]{RW}
}

\medskip

The two ideas mentioned above will be combined; at every local minimum of the
Brownian motion $ B $, a new random variable will appear from thin air. That
is, the binary extension, performed in Sect.~\ref{sec:6} at the global
minimum, will be performed at every local minimum, thus achieving locality and
stationarity required from a noise, while retaining the drift sensitivity
achieved in Sect.~\ref{sec:6} (as will be shown in Sect.~\ref{sec:8}).

A new random sign attached to a local minimum at $ \tau $ may be thought of
as a random choice of one of the two functions $ S : (\tau,\tau+\eps_1) \to
\{-1,+1\} $ constant on $ [\tau+\eps_{n+1}, \tau+\eps_n) $ (for each $ n $)
and such that
\begin{equation}\label{7.*}
S(\tau+\eps_n) = S(\tau+\eps_n-) f_n \( B(\tau+\eps_n) - B(\tau) \)
\end{equation}
(the numbers $ \eps_n \downarrow 0 $ and the functions $ f_n : \R \to
\{-1,+1\} $ being chosen appropriately). Given a time interval $ (0,t) $, for
each local minimizer $ \tau \in (0,t) $ we describe the new random
sign by the value $ S(t) $ (of the corresponding function $ S $), denoted
however by $ \eta_t (\tau) $. Relation \eqref{7.*} turns into the relation
\eqref{eq:2.5} between $ \eta_{s} (\tau) $ and $ \eta_{s+t} (\tau) $.

Before attaching something to the local minima we enumerate them. For every
time interval $ (a,b) \subset \R $ there exists a \emph{measurable
enumeration} of local minima on $ (a,b) $, --- a sequence of
\measurable{\F_{a,b}^\white} random variables $ \tau_1, \tau_2, \dots : \Om
\to (a,b) $ such that for almost all $ \om $ the Brownian path $ t \mapsto
B_t(\om) $ has a local minimum at each $ \tau_k(\om) $, no other local minima
exist on $ (a,b) $, and the numbers $ \tau_1(\om), \tau_2(\om), \dots $ are
pairwise different a.s. Here is a simple construction for $ (a,b) = (0,1) $
taken from \cite[2e]{Ts04}. First, $ \tau_1(\om) $ is the minimizer on the
whole $ (0,1) $ (unique a.s.). Second, if $ \tau_1(\om) \in (0,1/2) $ then $
\tau_2(\om) $ is the minimizer on $ (1/2,1) $, otherwise --- on $ (0,1/2)
$. Third, $ \tau_3(\om) $ is the minimizer on the first of the four intervals
$ (0,1/4) $, $ (1/4,1/2) $, $ (1/2,3/4) $ and $ (3/4,1) $ that contains
neither $ \tau_1(\om) $ nor $ \tau_2(\om) $. And so on.

All measurable enumerations $ (\tau'_k)_k $ result from one of them $
(\tau_k)_k $ in the sense that
\[
\tau'_k(\om) = \tau_{\si_\om(k)} (\om) \quad \text{a.s.}
\]
for some (unique, in fact) \emph{random permutation} $ \si : \Om \to S_\infty
$, that is, an \measurable{\F_{a,b}^\white} random variable $ \si $ valued in
the group $ S_\infty $ of all bijective maps $ \{1,2,\dots\} \to \{1,2,\dots\}
$ (equipped with its natural Borel \sif). See also \cite[2e]{Ts04}.

Each $ \tau_k $ is a measurable selector of the set of all local minimizers;
for short, let us say just a \emph{selected minimum.} Here is the general form
of a selected minimum $ \tau $ in terms of a given enumeration $ (\tau_k)_k $:
\begin{equation}
  \label{eq:2.*}
  \begin{aligned}
    & \tau(\om) = \tau_k(\om) \quad \text{for } \om \in A_k \, , \\
    & \quad \text{where } (A_1,A_2,\dots) \text{ is a countable measurable
      partition of $ \Om $.}
  \end{aligned}
\end{equation}
Every selected minimum may serve as (say) the first element of some
enumeration.

Given two adjacent time intervals $ (a,b) $ and $ (b,c) $, we may choose
measurable enumerations $ (\tau'_k)_k $ and $ (\tau''_k)_k $ of local minima
on $ (a,b) $ and $ (b,c) $ respectively, and combine them into a measurable
enumeration $ (\tau_k)_k $ on $ (a,c) $, say,
\begin{equation}
  \label{eq:2.1}
\tau_{2k-1} = \tau'_k \, , \quad \tau_{2k} = \tau''_k \quad \text{for }
  k=1,2,\dots
\end{equation}
taking into account that the point $ b $ is a.s.\ not a local minimizer.

Now we can attach independent random signs to the local
minima. Let $ \Om_1^\white \subset C[0,1] $ be a set of full Wiener
measure such that for every $ \om^\white_1 \in \Om_1^\white $ the set
$ \LocMin (\om^\white_1) $ of all local minimizers of the path $
\om^\white_1 $ is a dense countable subset of $ (0,1) $. We introduce the
set $ \{-1,+1\}^{ \LocMin (\om^\white_1) } $ of all functions $ \eta_1
: \LocMin (\om^\white_1) \to \{-1,+1\} $ and consider the disjoint union $
\Om_1 $ of these sets over all $ \om^\white_1 $,
\[
\Om_1 = \{ ( \om^\white_1, \eta_1 ) : \om^\white_1 \in
\Om_1^\white, \, \eta_1 \in \{-1,+1\}^{ \LocMin (\om^\white_1) } \} \, .
\]

Every measurable enumeration $ (\tau_k)_k $ of the local minima on $ (0,1) $
gives us a one-to-one correspondence
\begin{gather*}
\Om_1 \leftrightarrow \Om_1^\white \times \{-1,+1\}^\infty \, , \\
( \om^\white_1, \eta_1 ) \leftrightarrow \( \om^\white_1, \(
  \eta_1 ( \tau_k(\om^\white_1)) \)_k \) \, ;
\end{gather*}
here $ \{-1,+1\}^\infty = \{-1,+1\}^{\{1,2,\dots\}} $ is the set of all infinite sequences
of numbers $ \pm1 $. (As usual, a set of Wiener measure $ 0 $ in $
\Om_1^\white $ may be neglected.) We take the uniform probability
distribution $ m $ on $ \{-1,+1\} $
(giving equal probabilities $ 0.5 $ to $ -1 $ and $ +1 $), equip $
\{-1,+1\}^\infty $ with the product measure $ m^\infty $, and $ \Om_1^\white
\times \{-1,+1\}^\infty $ --- with the Wiener measure multiplied by $ m^\infty
$. Then,
using the one-to-one correspondence, we transfer the probability measure (and
the underlying \sif) to $ \Om_1 $. The choice of an enumeration $
(\tau_k)_k $ does not matter, since $ m^\infty $ is invariant under
permutations.

Now $ \Om_1 $ is a probability space. Similarly, $ \Om_t $ becomes a
probability space for every $ t \in (0,\infty) $. Given $ s,t \in (0,\infty)
$, we get a natural isomorphism
\begin{gather}
\Om_s \times \Om_t \longleftrightarrow \Om_{s+t} \, , \label{eq:2.2}
 \\
\( (\om_s^\white, \eta_s), (\om_t^\white, \eta_t) \)
 \longleftrightarrow (\om_{s+t}^\white, \eta_{s+t}) \notag
\end{gather}
where $ (\om_s^\white, \om_t^\white) \longleftrightarrow
\om_{s+t}^\white $ is the usual composition of Brownian paths, and 
\begin{equation}
  \label{eq:2.3}
\eta_{s+t} (\tau) = \begin{cases}
\eta_s (\tau) &\text{if } \tau < s , \\
\eta_t (\tau-s) &\text{if } \tau > s .
\end{cases}
\end{equation}
(The notation is not good for the case $ s=t $, since $ \om_s $ and $ \om_t $
are still treated as different variables; hopefully it is not too confusing.)
The
composition $ (\eta_s,\eta_t) \longleftrightarrow \eta_{s+t} $ is
described conveniently in terms of an enumeration of the form \eqref{eq:2.1}
for $ a=0 $, $ b=s $, $ c=s+t $:
\begin{equation}
  \eta_{s+t} (\tau_{2k-1}) = \eta_s (\tau'_k) \, , \quad \eta_{s+t}
  (\tau_{2k}) = \eta_t (\tau''_k-s) \quad \text{for } k = 1,2,\dots
\end{equation}
(of course, all these $ \eta $ and $ \tau $ depend implicitly on the
underlying $ \om^\white $).

We have a noise
(an extension of the white noise). It is described above via probability
spaces $ \Om_t $ satisfying $ \Om_s \times \Om_t = \Om_{s+t} $ rather than
sub-\sif s $ \F_{s,t} $ (on a single $ \Om $) satisfying $ \F_{r,s} \otimes
\F_{s,t} = \F_{r,t} $, but these are two equivalent languages (see \cite[3c1
and 3c6]{Ts04}), and the corresponding Arveson system is just $ H_t =
L_2(\Om_t) $.

However, it is not yet the new, drift sensitive noise that we need. Rather, it
is Warren's noise of splitting. The binary extension performed at each $ \tau
$ should follow the construction of Sect.~\ref{sec:6}.
To this end we retain the probability spaces $ \Om_t $ constructed before,
but replace the straightforward isomorphisms \eqref{eq:2.2}--\eqref{eq:2.3}
with less evident, `twisted' isomorphisms. Namely, \eqref{eq:2.3} is replaced
with
\begin{gather}
 \eta_{s+t} (\tau) = \eta_{t} (\tau-s) \quad \text{if } \tau > s \, ,
  \label{eq:2.4} \\
 \eta_{s+t} (\tau) = \eta_s (\tau) \prod_{n:\tau+\eps_n\in(s,s+t]} f_n \(
  B(\tau+\eps_n) - B(\tau) \) \quad \text{if } \tau < s \, . \label{eq:2.5}
\end{gather}
As before, all these $ \eta $ and $ \tau $ depend implicitly on the underlying
$ \om^\white $, and $ B(s)(\om_t^\white) = \om_t^\white(s) $ for $ s \in [0,t]
$.

The new noise is thus constructed. Its parameters $ (\eps_n)_n $ and $ (f_n)_n
$ will be chosen later. (In fact, $ \eps_n = 2 \cdot 3^{-n-1} $, and $ f_n $
are given by Lemma \ref{5.5new}.)

The classical part of the new noise is exhausted by the white noise, which can
be proved via the predictable representation property, see \cite[4d]{Ts04}.

In order to examine the impact of drift on the new noise we need the relation
\begin{equation}\label{drift}
\LocMin (\om_t^\white) = \LocMin \( \theta_t^\la(\om_t^\white) \)
\end{equation}
(for all $ t,\la $ and almost all $ \om_t^\white \in \Om_t^\white $); as
before, $ \theta_t^\la : C[0,t] \to C[0,t] $ is the drift transformation,
$ (\theta_t^\la b)(s) = b(s)-2\la s $. The relation \eqref{drift} follows from
the well-known fact that all local minima of the Brownian motion are sharp
(a.s.) in the sense that
\[
\frac{ B(t) - B(\tau) }{ |t-\tau| } \to \infty \quad \text{as } t \to \tau, \,
t \ne \tau
\]
whenever $ \tau $ is a local minimizer. See \cite[Sect.~2.10, Items
7,8]{IM}. (In fact, $ |t-\tau| $ may be replaced with $ \sqrt{|t-\tau|} /
\ln^2 |t-\tau| $.)

It is easy to guess that a drift corresponds to a shift of the functions $ f_n
$. The proof (rather boring) is given below.

\begin{lemma}\label{10.10}
Let numbers $ \la \in \R $, $ \eps_n \downarrow 0 $ and Borel functions $ f_n,
g_n : \R \to \{-1,+1\} $ satisfy
\[
g_n(x) = f_n(x+2\la\eps_n) \quad \text{for all $ x \in \R $ and $ n $} \, .
\]
Let two extensions of the white noise be constructed as before, one
corresponding to $ (f_n)_n $ and $ (\eps_n)_n $, the other corresponding to $
(g_n)_n $ and $ (\eps_n)_n $. Then the second extension results from the first
one by the drift $ 2\la $ (as defined in Sect.~\ref{sec:1-5}), up to
isomorphism of extensions.
\end{lemma}

\begin{proof}
The probability spaces $ \Om_t $ and measure preserving maps $ \Om_t \to
\Om_t^\white $ are the same for both extensions, however, the corresponding
isomorphisms $ \al_f, \al_g : \Om_s \times \Om_t \to \Om_{s+t} $ differ; $
\al_f $, used in the first extension, involves $ f_n $ (recall
\eqref{eq:2.5}), while $ \al_g $, used in the second extension, involves $ g_n
$ instead of $ f_n $.

We introduce the third extension, resulting from the first one by the drift $
2\la $, and seek an isomorphism between the second and third extensions.

The third extension uses the same $ \Om_t $ but with probability measures $
P'_t $ different from the probability measures $ P_t $ used by the first
and second extensions; namely,
\[
\frac{ \D P'_t }{ \D P_t } = D_t = \exp ( 2\la B_t - 2 \la^2 t ) \, .
\]
The white noise extended by the third extension is generated by the Brownian
motion $ B'_t = B_t - 2\la t $. Note also that the third extension uses $
\al_f $.

The probability space $ \Om_t $ consists of pairs $ (\om_t^\white, \eta_t) $
where $ \om_t^\white \in \Om_t^\white \subset C[0,t] $ and $ \eta_t \in
\{-1,+1\}^{\LocMin(\om_t^\white)} $. The drift transformation $ \theta_t^\la $
may be treated as a measure preserving map
\[
\theta_t^\la : ( \Om_t^\white, D_t \cdot \W_t ) \to ( \Om_t^\white, \W_t ) \,
.
\]
Using \eqref{drift} we define $ \ti\theta_t^\la : \Om_t \to \Om_t $ by $
\ti\theta_t^\la (\om_t^\white, \eta_t) = (\theta_t^\la \om_t^\white, \eta_t) $
and get a measure preserving map
\[
\ti\theta_t^\la : (\Om_t,P'_t) \to (\Om_t,P_t) \, .
\]
Clearly, $ B'_s = B_s \circ \ti\theta_t^\la $ for $ s \in [0,t] $. It remains
to check that $ \ti\theta_s^\la \times \ti\theta_t^\la = \ti\theta_{s+t}^\la $
in the sense that the diagram
\[
\xymatrix@C+1cm{
 \Om_s \times \Om_t \ar[d]^{\al_f}
  \ar[r]^{\ti\theta_s^\la \times \ti\theta_t^\la} &
 \Om_s \times \Om_t \ar[d]^{\al_g}
\\
 \Om_{s+t} \ar[r]^{\ti\theta_{s+t}^\la} & \Om_{s+t}
}
\]
is commutative. Let $ \om_s = (\om_s^\white,\eta_s) \in \Om_s $ and $ \om_t =
(\om_t^\white,\eta_t) \in \Om_t $. We have $ \al_f (\om_s,\om_t) =
(\om_{s+t}^\white, \eta_{s+t}) $, where $ \om_{s+t}^\white $ is the usual
composition of $ \om_s^\white $ and $ \om_t^\white $, while $ \eta_{s+t} $ is
obtained from $ \eta_s $ and $ \eta_t $ according to \eqref{eq:2.4},
\eqref{eq:2.5}. Thus,
\[
\ti\theta_{s+t}^\la \( \al_f (\om_s,\om_t) \) = \ti\theta_{s+t}^\la (
\om_{s+t}^\white, \eta_{s+t} ) = \( \theta_{s+t}^\la (\om_{s+t}^\white),
\eta_{s+t} \) \, .
\]
On the other hand,
\[
( \ti\theta_s^\la \times \ti\theta_t^\la ) ( \om_s, \om_t ) = \(
\ti\theta_s^\la (\om_s), \ti\theta_t^\la (\om_t) \) = \( ( \theta_s^\la
(\om_s^\white), \eta_s ), ( \theta_t^\la (\om_t^\white), \eta_t ) \) \, .
\]
Clearly, $ \al_g \( ( \theta_s^\la (\om_s^\white), \eta_s ), ( \theta_t^\la
(\om_t^\white), \eta_t ) \) = \( \theta_{s+t}^\la (\om_{s+t}^\white),
\eta'_{s+t} \) $ for some $ \eta'_{s+t} $ (since $ \theta_s^\la \times
\theta_t^\la = \theta_{s+t}^\la $). Finally, $ \eta'_{s+t} = \eta_{s+t} $ by
\eqref{eq:2.4}, \eqref{eq:2.5} and the equality
\begin{multline*}
g_n \( \theta_{s+t}^\la (\om_{s+t}^\white) (\tau+\eps_n) - \theta_{s+t}^\la
 (\om_{s+t}^\white) (\tau) \) = \\
= g_n \( \om_{s+t}^\white (\tau+\eps_n) - \om_{s+t}^\white (\tau) - 2\la\eps_n
 \) = f_n \( \om_{s+t}^\white (\tau+\eps_n) - \om_{s+t}^\white (\tau) \) \, .
\end{multline*}
\end{proof}

\section[The binary extension inside the new noise]
 {\raggedright The binary extension inside the new noise}
\label{sec:8}
According to Sect.~\ref{sec:6}, the Brownian motion $ B $ leads to an
inductive system of probability spaces formed by the restrictions of $ B $ to
the time intervals $ [0,\tau+3^{-n}] \cap [0,1] $, where $ \tau $ is the
(global) minimizer of $ B $ on $ [0,1] $. Further, every sequence $ (f_n)_n $
of Borel functions $ f_n : \R \to \{-1,+1\} $ leads to a binary extension of
this inductive system. The extension is formed by the restrictions of $ B $
and $ S_f $ to $ [0,\tau+3^{-n}] \cap [0,1] $; here $ S_f : (\tau,1) \to
\{-1,+1\} $ is a random function satisfying \eqref{6.9}.

On the other hand, according to Sect.~\ref{sec:7}, $ (f_n)_n $ (in combination
with $ \eps_n = 2 \cdot 3^{-n-1} $) leads to a noise that extends the white
noise. The noise is formed by the Brownian motion $ B $ and the random
variables $ \eta_t(\tau) $; here $ \tau $ runs over all local minimizers of $
B $ on $ (0,t) $. In turn, the noise leads to an Arveson system that extends
the type $ I_1 $ Arveson system of the white noise.

These constructions of Sections \ref{sec:6} and \ref{sec:7} are related as
follows.

\begin{proposition}\label{10.1}
If two sequences $ (f_n)_n $, $ (g_n)_n $ of Borel functions $ \R \to
\{-1,+1\} $ lead to isomorphic extensions of the type $ I_1 $ Arveson system
(of the white noise), then they lead to isomorphic binary extensions of the
inductive system of probability spaces.
\end{proposition}

The proof is given after the proof of Prop.~\ref{10.4}.

\begin{proof}[Proof of Theorem \ref{1.8}]
The binary extension, constructed in Sect.~\ref{sec:6} using the functions $
f_n $ given by Lemma \ref{5.5new} (combined with Prop.~\ref{prop8}), is not
isomorphic to the extension that corresponds to the shifted functions $ g_n(x)
= f_n(x+3^{-n}c) $, unless $ c = 0 $. By Prop.~\ref{10.1}, $ (f_n)_n $ and $
(g_n)_n $ lead to nonisomorphic extensions (constructed in Sect.~\ref{sec:7})
of the type $ I_1 $ Arveson system (of the white noise). By Lemma \ref{10.10},
these nonisomorphic extensions result from one another by a drift. This drift
sensitivity implies Theorem \ref{1.8} by Corollary \ref{cor2}.
\end{proof}

Comparing \eqref{6.9} and \eqref{eq:2.5} we see that the function $ t \mapsto
\eta_t (\tau) $ behaves like the function $ S_f $. (Here $ \tau $ is the
global minimizer of $ B $ on $ [0,1] $.) In other words, we may let for some
(therefore, all) $ n $ such that $ \tau + 3^{-n} < 1 $,
\[
S_f (\tau+3^{-n}) = \eta_{\tau+3^{-n}} (\tau) \, ,
\]
thus defining a measure preserving map from the probability space $
\Om^{\noise(f)} $ of the noise on $ [0,1] $ to the probability space $
\Om^{\bin(f)} $ of the binary extension on $ [0,1] $;
\[
\Om^{\noise(f)} \to \Om^{\bin(f)} \, ;
\]
here $ f = (f_n)_n $ is the given sequence of functions. Accordingly, we have
a natural embedding of Hilbert spaces,
\[
L_2 ( \Om^{\bin(f)} ) \to L_2 ( \Om^{\noise(f)} ) \, .
\]
Striving to prove Prop.~\ref{10.1} we assume existence of an isomorphism $
\Theta = (\Theta_t)_t $ between the two Arveson systems,
\begin{gather}
\Theta_t : L_2 ( \Om_t^{\noise(f)} ) \to L_2 ( \Om_t^{\noise(g)} ) \,
 , \label{10.2} \\
\Theta_t \text{ is trivial on } L_2 ( \Om_t^{\white} ) \, . \label{10.3}
\end{gather}
Note that $ \Om^{\noise(f)} = \Om_1^{\noise(f)} $.

\begin{proposition}\label{10.4}
$ \Theta_1 $ maps the subspace $ L_2 ( \Om^{\bin(f)} ) $ of $ L_2 (
\Om^{\noise(f)} ) $ onto the subspace $ L_2 ( \Om^{\bin(g)} ) $ of $ L_2 (
\Om^{\noise(g)} ) $.
\end{proposition}

The proof is given after Lemma \ref{10.5}.

The structure of $ L_2(\Om^{\noise(f)}) $ is easy to describe:
\[
L_2(\Om^{\noise(f)}) = H_0^f \oplus H_1^f \oplus H_2^f \oplus \dots \, ,
\]
where $ H_n^f $ (called the $n$-th superchaos space) consists of the random
variables of the form
\begin{gather*}
\sum_{k_1<\dots<k_n} \eta_1(\tau_{k_1}) \dots \eta_1(\tau_{k_n})
 \phi_{k_1,\dots,k_n} \, , \\
\phi_{k_1,\dots,k_n} \in L_2(\Om_1^\white) \, , \quad \sum_{k_1<\dots<k_n} \|
\phi_{k_1,\dots,k_n} \|^2 < \infty \, ,
\end{gather*}
where $ (\tau_k)_k $ is a measurable enumeration of the local minimizers of $ B
$ on $ (0,1) $ (the choice of the enumeration does not matter). See
\cite[(3.1)]{Ts06} for the case $ f_n(\cdot) = 1 $ (Warren's noise of
splitting); the same argument works in general. Note that $ H_0^f =
L_2(\Om_1^\white) $.

It is well-known that the superchaos spaces may be described in terms of the
Arveson system, and therefore $ \Theta_1 $ maps $ H_n^f $ onto $ H_n^g $. We
need the first superchaos space only; here is a simple argument for this case:
\[
H_1 = \{ \psi \in L_2(\Om^{\noise}) : \forall t \in (0,1) \;\> \psi = Q_{0,t}
\psi + Q_{t,1} \psi \} \, ;
\]
here $ Q_{0,t} $ is the orthogonal projection of the space $ L_2(\Om^{\noise})
= L_2(\Om_{t}^{\noise}) \otimes L_2(\Om_{1-t}^{\noise}) $ onto the subspace
$ L_2(\Om_{t}^{\noise}) \otimes L_2(\Om_{1-t}^{\white}) $, and $ Q_{t,1} $
--- onto $ L_2(\Om_{t}^{\white})\linebreak[0]
\otimes L_2(\Om_{1-t}^{\noise}) $. We have
\begin{multline*}
\Theta_1 \( L_2(\Om_{t}^{\noise(f)}) \otimes L_2(\Om_{1-t}^{\white}) \) = \\
= \Theta_t \( L_2(\Om_{t}^{\noise(f)}) \) \otimes \Theta_{1-t} \(
 L_2(\Om_{1-t}^{\white}) \) = \\
= L_2(\Om_{t}^{\noise(g)}) \otimes L_2(\Om_{1-t}^{\white})
\end{multline*}
by \eqref{10.3}; therefore $ \Theta_1 Q_{0,t}^f = Q_{0,t}^g \Theta_1
$. Similarly, $ \Theta_1 Q_{t,1}^f = Q_{t,1}^g \Theta_1 $. It follows that
\[
\Theta_1 H_1^f = H_1^g \, .
\]

Similarly, $ L_2(\Om_t^{\noise}) = H_0(t) \oplus H_1(t) \oplus H_2(t) \oplus
\dots $ (the upper index, be it $ f $ or $ g $, is omitted). Identifying $
L_2(\Om_1) $ with $ L_2(\Om_t) \otimes L_2(\Om_{1-t}) $ we have
\[
H_1 = \underbrace{ H_1(t) \otimes H_0(1-t) }_{Q_{0,t} H_1} \oplus \underbrace{
H_0(t) \otimes H_1(1-t) }_{Q_{t,1} H_1} \, .
\]

The commutative algebra $ L_\infty (\Om_1^\white) $ acts naturally on $
H_1 $:
\[
h \cdot \sum_k \eta_1(\tau_k) \phi_k = \sum_k \eta_1(\tau_k) h \cdot \phi_k
\quad \text{for } h \in L_\infty (\Om_1^\white) \, .
\]
Also the commutative algebra $ L_\infty (0,1) $ acts naturally on $ H_1 $. In
particular, $ \One_{(0,t)} $ acts as $ Q_{0,t} $, and $ \One_{(t,1)} $ acts as
$ Q_{t,1} $. In general,
\[
h \cdot \sum_k \eta_1(\tau_k) \phi_k = \sum_k \eta_1(\tau_k) h(\tau_k) \phi_k
\quad \text{for } h \in L_\infty (0,1) \, .
\]
(The choice of enumeration $ (\tau_k)_k $ does not matter.) The two actions
commute, and may be combined into the action of $ L_\infty(\mu) $ (on $ H_1
$) for some measure $ \mu $ on $ \Om_1^\white \times (0,1) $:
\[
h \cdot \sum_k \eta_1(\tau_k) \phi_k(\cdot) = \sum_k \eta_1(\tau_k)
h(\cdot,\tau_k) \phi_k(\cdot) \quad \text{for } h \in L_\infty (\mu) \, .
\]
The measure $ \mu $ may be chosen as
\[
\int h \, \D\mu = \Ex \sum_k \frac1{k^2} h(B,\tau_k)
\]
(or anything equivalent).

\begin{lemma}\label{10.5}
The diagram
\[
\xymatrix@C+1cm{
 H_1^f \ar[d]^{h} \ar[r]^{\Theta_1} & H_1^g \ar[d]^{h}
  \\
 H_1^f \ar[r]^{\Theta_1} & H_1^g
}
\]
is commutative for every $ h \in L_\infty(\mu) $.
\end{lemma}

\begin{proof}
Given $ [a,b] \subset [0,1] $, we define a subalgebra $ \Ga(a,b) \subset
L_\infty(\mu) $ as consisting of the functions of the form
\[
h(\om_{0,1}^\white,t) = \begin{cases}
 h'(\om_{0,a}^\white,\om_{b,1}^\white) &\text{for } t \in (a,b), \\
 0 &\text{for } t \in (0,a) \cup (b,1),
\end{cases}
\]
where $ h' \in L_\infty ( \Om_a^\white \times \Om_{1-b}^\white ) $, and 
$ \om_{0,1}^\white \in \Om_1^\white $ is treated as the triple $ (
\om_{0,a}^\white, \om_{a,b}^\white, \om_{b,1}^\white ) $ according to the
natural isomorphism between $ \Om_1^\white $ and $ \Om_a^\white \times
\Om_{b-a}^\white \times \Om_{1-b}^\white $. For each $ n=1,2,\dots $ we define
a subalgebra $ \Ga_n \subset L_\infty(\mu) $ by
\[
\Ga_n = \sum_{k=1}^{2^n} \Ga \Big( \frac{k-1}{2^n}, \frac{k}{2^n} \Big) \, .
\]
It is easy to see that $ \Ga_n $ corresponds to a measurable partition; in
other words, $ \Ga_n = L_\infty ( \Om_1^\white \times (0,1), \Ec_n, \mu )
$ for some sub-\sif\ $ \Ec_n $ of the \sif\ $ \Ec $ of all \measurable{\mu}
sets. We have $ \Ec_n \uparrow \Ec $, that is, $ \Ec_1 \subset \Ec_2 \subset
\dots $ and $ \Ec $ is the least sub-\sif\ containing all $ \Ec_n $, which
follows from the fact that $ \Ga_1 \cup \Ga_2 \cup \dots $ contains a
countable set that separates points of $ \Om_1^\white \times (0,1) $.

If $ \Theta_1 h_n = h_n \Theta_1 $ (as operators $ H_1^f \to H_1^g $) for
all $ n $, and $ h_n \to h $ almost everywhere, and $ \sup_n \| h_n
\|_\infty < \infty $, then $ \Theta_1 h = h \Theta_1 $. Thus, it is
sufficient to prove the equality $ \Theta_1 h = h \Theta_1 $ for all $ h
\in \Ga_1 \cup \Ga_2 \cup \dots $. Without loss of generality we may assume
that $ h \in \Ga(a,b) $ for some $ a,b $. Moreover, I assume that $ b=1 $,
leaving the general case to the reader. Thus, $ h(\om_{0,1}^\white,t) =
h'(\om_{0,a}^\white) \One_{(a,1)} (t) $.

We recall that $ \Theta_1 = \Theta_a \otimes \Theta_{1-a} $, $ H_1 = H_1(a)
\otimes H_0(1-a) \oplus H_0(a) \otimes H_1(1-a) $, and note that $ \Theta_1 \(
H_1^f(a) \otimes H_0^f(1-a) \) = H_1^g(a) \otimes H_0^g(1-a) $, $ \Theta_1 \(
H_0^f(a) \otimes H_1^f(1-a) \) = H_0^g(a) \otimes H_1^g(1-a) $. The subspaces
$ H_1^f(a) \otimes H_0^f(1-a) $ and $ H_1^g(a) \otimes H_0^g(1-a) $ are
annihilated by $ h $, thus, $ \Theta_1 h $ and $ h \Theta_1 $ both vanish on
$ H_1^f(a) \otimes H_0^f(1-a) $. On the other subspace, $ H_0^f(a) \otimes
H_1^f(1-a) $, $ h $ acts as $ h' \otimes \One $, while $ \Theta_1 $ acts as $
\One \otimes \Theta_{1-a} $. Therefore $ \Theta_1 h = h \Theta_1 $.
\end{proof}

\begin{proof}[Proof of Prop.~\ref{10.4}]
We apply Lemma \ref{10.5} to the function $ h \in L_\infty(\mu) $ defined by
\[
h (\om_1^\white,t) = \begin{cases}
 1 & \text{if } \tau(\om_1^\white)=t, \\
 0 & \text{otherwise},
\end{cases}
\]
where $ \tau $ is the (global) minimizer on $ (0,1) $. This function acts on $
H_1 $ as the projection onto the subspace $ \{ \eta_1(\tau) \phi : \phi \in
L_2(\Om_1^\white) \} = L_2 (\Om^\bin) \ominus L_2 (\Om^\white) $.
\end{proof}

\begin{proof}[Proof of Prop.~\ref{10.1}]
We have two binary extensions, $ (\Om_n^{\bin(f)}, \ti\be_n^f)_n $ and $
(\Om_n^{\bin(g)},\linebreak[0]
\ti\be_n^g)_n $, of an inductive system $ (\Om_n^\white,
\be_n)_n $ of probability spaces (according to $ (\ga_n^f)_n $ and $
(\ga_n^g)_n $ respectively). Their isomorphism is ensured by Lemma \ref{6.22},
provided that Condition \ref{6.22}(b) is satisfied by some unitary operators 
$ \Theta_n^\bin : L_2(\Om_n^{\bin(f)}) \to L_2(\Om_n^{\bin(g)}) $. Using the
natural embeddings $ L_2(\Om_n^{\bin(f)}) \subset L_2(\Om^{\bin(f)}) $ and $
L_2(\Om_n^{\bin(g)}) \subset L_2(\Om^{\bin(g)}) $ we define all $
\Theta_n^\bin $ as restrictions of a single operator $ \Theta^\bin :
L_2(\Om^{\bin(f)}) \to L_2(\Om^{\bin(g)}) $. Using Prop.~\ref{10.4} we define
$ \Theta^\bin $ as the restriction of $ \Theta_1 $ to $ L_2(\Om^{\bin(f)})
$. It remains to prove that
\begin{gather}
\Theta^\bin \( L_2(\Om_n^{\bin(f)}) \) = L_2(\Om_n^{\bin(g)}) \, ,
 \label{three-a} \\
\begin{gathered}
 \Theta_n^\bin \text{ intertwines the actions of }
 L_\infty(\Om_n^\white)\qquad\qquad \\
 \qquad\qquad\qquad\qquad \text{ on } L_2(\Om_n^{\bin(f)}) \text{ and }
 L_2(\Om_n^{\bin(g)}) \, ,
\end{gathered}
 \label{three-b} \\
\Theta_n^\bin \text{ is trivial on } L_2(\Om_n^\white)  \label{three-c}
\end{gather}
for all $ n $.

By \eqref{10.3}, $ \Theta^\bin $ is trivial on $ L_2(\Om^\white) $;
\eqref{three-c} follows.

By Lemma \ref{10.5}, $ \Theta^\bin $ intertwines the actions of $
L_\infty(\Om^\white) $ on $ L_2(\Om^{\bin(f)}) $ and $ L_2(\Om^{\bin(g)}) $;
\eqref{three-b} follows.

The proof of \eqref{three-a} is the point of Prop.~\ref{10.6} below.
\end{proof}

\begin{proposition}\label{10.6}
The operator $ \Theta^\bin $ maps the subspace $ L_2(\Om_n^{\bin(f)}) \subset
L_2(\Om^{\bin(f)}) $ onto the subspace $ L_2(\Om_n^{\bin(g)}) \subset
L_2(\Om^{\bin(g)}) $.
\end{proposition}

The proof is given after Lemma \ref{11.7}.

Recall that the elements of $ L_2(\Om_n^{\bin(f)}) $ are functions of the
restrictions of $ B $ and $ S_f $ to $ [0,\tau+3^{-n}] \cap [0,1] $.

For a given $ t \in (0,1) $ we consider the sub-\sif\ $ \F_t^f $ on $
\Om^{\bin(f)} $, generated by the restrictions of $ B $ and $ S_f $ to $ [0,t]
$. The elements of the subspace $ L_2(\Om^{\bin(f)}, \F_t^f) $ are functions
of $ B|_{[0,t]} $ and $ S_f|_{[0,t]} $.

We know that $ L_\infty (\Om^\white) $ acts on $ L_2 (\Om^{\bin(f)}) $. In
particular, for $ 0<r<s<1 $, the function $ \One_{(r,s)} (\tau) $ (that is,
the indicator of $ \{ \om^\white : r<\tau(\om^\white)<s \} $) acts as the
projection onto a subspace $ H_{r,s}^f $ of $ L_2 (\Om^{\bin(f)}) $. The same
holds for $ g $. We have $ \Theta^\bin (H_{r,s}^f) \subset H_{r,s}^g $, since
$ \Theta^\bin $ intertwines the two actions of $ L_\infty(\Om^\white) $. We
define
\[
H_{r,s,t}^f = H_{r,s}^f \cap L_2 ( \Om^{\bin(f)}, \F_t^f ) \quad \text{for }
0<r<s<t<1 \, .
\]

\begin{lemma}\label{11.7}
$ \Theta^\bin (H_{r,s,t}^f) \subset H_{r,s,t}^g $.
\end{lemma}

\begin{proof}
The binary extension $ \Om^{\bin(f)} $ is constructed on the time interval $
(0,1) $, but the same can be made on the time interval $ (0,t) $, giving a
binary extension $ \Om^{\bin(f,t)} $ of $ \Om_t^\white $, using the (global)
minimizer $ \tau_t $ on $ (0,t) $; sometimes $ \tau_t = \tau_1 $, sometimes $
\tau_t \ne \tau_1 $.

The binary extension $ \Om^{\bin(f)} $ is the product (recall
Def.~\ref{4-8.1}) of two binary extensions, $ \Om^{\bin(f,t)} $ and $
\Om^{\bin(f,1-t)} $, according to the set $ A \subset \Om_1^\white =
\Om_t^\white \times \Om_{1-t}^\white $,
\[
A = \{ \om_1^\white : \tau_1(\om_1^\white) = \tau_t(\om_1^\white) \} \, .
\]
We know that $ \Theta_1 = \Theta_t \otimes \Theta_{1-t} $. Similarly to
Prop.~\ref{10.4}, $ \Theta_t ( L_2(\Om^{\bin(f,t)}) ) = L_2(\Om^{\bin(g,t)})
$; we define $ \Theta^{\bin,t} : L_2(\Om^{\bin(f,t)}) \to L_2(\Om^{\bin(g,t)})
$ as the restriction of $ \Theta_t $ and observe that $ \Theta^\bin $ is the
restriction of $ \Theta^{\bin,t} \otimes \Theta^{\bin,1-t} $ to $ L_2
(\Om^{\bin(f)}) \subset L_2(\Om^{\bin(f,t)}) \otimes L_2(\Om^{\bin(f,1-t)}) $
(recall \eqref{4-8.*}).

By Lemma \ref{4-8.2}, $ \Theta^\bin (L_2(\ti A,\F_1)) = L_2(\ti A',\F'_1) $,
where the sets $ \ti A \subset \Om^{\bin(f)} $, $ \ti A' \subset \Om^{\bin(g)}
$ correspond to the inequality $ \tau < t $, the sub-\sif\ $ \F_1 $ on $ \ti A
$ is induced by the sub-\sif\ $ \F_t^f $ on $ \Om^{\bin(f)} $, and $ \F'_1 $
on $ \ti A' $ --- by $ \F_t^g $.

Taking into account that $ (r,s) \subset (0,t) $ we get $ H_{r,s,t} \subset
L_2(\ti A,\F_1) $. Therefore $ \Theta^\bin (H^f_{r,s,t}) \subset L_2(\ti
A',\F'_1) $. On the other hand, $ \Theta^\bin (H^f_{r,s,t}) \subset
\Theta^\bin (H^f_{r,s}) \subset H^g_{r,s} $. It remains to note that $ L_2(\ti
A',\F'_1) \cap H^g_{r,s} \subset H^g_{r,s,t} $.
\end{proof}

\begin{proof}[Proof of Prop.~\ref{10.6}]
If $ r,s,t $ and $ n $ satisfy $ t \le r + 3^{-n} $ then $ H_{r,s,t}^g \subset
L_2 (\Om_n^{\bin(g)}) $ (since $ t \le \tau(\cdot) + 3^{-n} $ for all relevant
points), and therefore $ \Theta^\bin (H_{r,s,t}^f) \subset L_2
(\Om_n^{\bin(g)}) $.

The elements of $ L_2(\Om_n^{\bin(f)}) $ are functions of the restrictions of
$ B $ and $ S_f $ to $ [0,\tau+3^{-n}] \cap [0,1] $. For every $ N $ such that
$ \frac1N < 3^{-n} $ consider the functions of the restrictions of $ B $ and $
S_f $ to $ [0,\tau+3^{-n}-\frac1N] \cap [0,1] $; these are $
L_2(\Om_n^{\bin(f)},\Ec_N) $ for some sub-\sif\ $ \Ec_N $, and
\[
\bigcup_N L_2(\Om_n^{\bin(f)},\Ec_N) \quad \text{is dense in }
L_2(\Om_n^{\bin(f)}) \, ,
\]
since $ \Ec_N \uparrow \Ec $ (a similar argument is used in the proof of Lemma
\ref{10.5}; note that $ S_f $ jumps at $ \tau + 2 \cdot 3^{-n} $, not $ \tau +
3^{-n} $). In order to prove Prop.\ \ref{10.6} it remains to prove that
\[
\Theta^\bin \( L_2(\Om_n^{\bin(f)},\Ec_N) \) \subset L_2(\Om^{\bin(g)})
\]
for all $ N $ (satisfying $ \frac1N < 3^{-n} $).

Clearly,
\[
L_2(\Om^{\bin(f)}) = H_{0,\frac1N} \oplus \dots \oplus H_{\frac{N-1}N,1}
\]
(for every $ N $). Every $ \psi \in L_2(\Om^{\bin(f)}) $ is of the form
\[
\psi = \psi_1 + \dots + \psi_N \, , \quad \psi_k \in H_{\frac{k-1}N,\frac k N}
\, .
\]
If $ \psi \in L_2(\Om_n^{\bin(f)},\Ec_N) $ then $ \psi_k \in L_2
(\Om^{\bin(f)}, \F_{\frac{k-1}N+3^{-n}}) $ (since $ \tau(\cdot) + 3^{-n} -
\frac1N < \frac{k-1}N + 3^{-n} $ for all relevant points), thus, $ \psi_k \in
H^f_{ \frac{k-1}N, \frac k N, \frac{k-1}N+3^{-n} } $. Taking into account that
$ \Theta^\bin \( H^f_{ \frac{k-1}N, \frac k N, \frac{k-1}N+3^{-n} } \) \subset
L_2 (\Om_n^{\bin(g)}) $ we see that $ \Theta^\bin (\psi) \in L_2
(\Om_n^{\bin(g)}) $.
\end{proof}

\bigskip
\filbreak
{
\small
\begin{sc}
\parindent=0pt\baselineskip=12pt
\parbox{4in}{
Boris Tsirelson\\
School of Mathematics\\
Tel Aviv University\\
Tel Aviv 69978, Israel
\smallskip
\par\quad\href{mailto:tsirel@post.tau.ac.il}{\tt
 mailto:tsirel@post.tau.ac.il}
\par\quad\href{http://www.tau.ac.il/~tsirel/}{\tt
 http://www.tau.ac.il/\textasciitilde tsirel/}
}

\end{sc}
}
\filbreak

\end{document}